\documentclass[11pt,fleqn]{article}
\usepackage{amssymb}
\usepackage{amsfonts}
\usepackage{amsmath}
\usepackage{amsthm}
\usepackage{graphicx}
\usepackage{epsfig}
\usepackage{psfrag}
\usepackage{dsfont}
\usepackage{amstext, latexsym, color}
\makeatletter
\xdef\@endgadget#1{{\unskip\nobreak\hfil\penalty50\hskip1em\hbox{}\nobreak
   \hfil#1\parfillskip=0pt\finalhyphendemerits=0\par}}
\def\@qedsymbol{${}_\blacksquare$}
\def\qed{\@endgadget{\@qedsymbol}}
\newtheorem{lemma}{Lemma}[section]
\newtheorem{theorem}[lemma]{Theorem}
\newtheorem{corollary}[lemma]{Corollary}

\newtheorem{definition}[lemma]{Definition}

\newtheorem{proposition}[lemma]{Proposition}
\newtheorem{remark}[lemma]{Remark}
\newtheorem{assumption}[lemma]{Assumption}
\newcommand{\mR}{\mathbb{R}}

\newcommand{\pperp}{\perp \!\!\!\perp}
\DeclareMathOperator{\im}{im}

\DeclareMathOperator{\spa}{span}

\def\BibTeX{{\rm B\kern-.05em{\sc i\kern-.025em b}\kern-.08em
   T\kern-.1667em\lower.7ex\hbox{E}\kern-.125emX}}

\title{\LARGE \bf
Port-Hamiltonian systems on graphs}

\author{A.J. van der Schaft and B.M. Maschke
\thanks{A.J. van der Schaft is with the Johann Bernoulli Institute for Mathematics and Computer
Science, University of Groningen, PO Box 407, 9700 AK, the
Netherlands
       {\tt\small A.J.van.der.Schaft@rug.nl}}
\thanks{B.M. Maschke is with the Laboratoire d'Automatique et de Genie des Proc\'ed\'es, Universit\'e Claude Bernard Lyon-1, F-69622 Villeurbanne, Cedex, France
       {\tt\small maschke@lagep.univ-lyon1.fr}}
\\
\\
August 25, 2012\footnote{Revised version of arXiv:1107.2006, July 12, 2011.}}

\date{}
\begin{document}

\maketitle
\thispagestyle{empty}
\pagestyle{empty}

\begin{abstract}     
In this paper we present a unifying geometric and compositional framework for modeling complex physical network dynamics as port-Hamiltonian systems on open graphs. Basic idea is to associate with the incidence matrix of the graph a Dirac structure relating the flow and effort variables associated to the edges, internal vertices, as well as boundary vertices of the graph, and to formulate energy-storing or energy-dissipating relations between the flow and effort variables of the edges and internal vertices. This allows for state variables associated to the edges, and formalizes the interconnection of networks. Examples from different origins such as consensus algorithms are shown to share the same structure. It is shown how the identified Hamiltonian structure offers systematic tools for the analysis of the resulting dynamics.
\end{abstract}

\section{Introduction}
Discrete topological structures arise abundantly in physical
systems modeling. A classical approach to the analysis of electrical circuits,
dating back to Kirchhoff, is based on the circuit graph. Similar approaches
apply to many other cases, including e.g. mass-spring-damper mechanical systems, multi-body
systems, hydraulic networks, chemical reaction networks, and power systems. 
Common feature is that the discrete structures, in particular graphs, are blended with dynamical relations, leading to various
sorts of \textit{network dynamics}. 

During the last two decades network dynamics has received ever-increasing attention, with inputs,
among others, from graph theory, multi-agent
systems, dynamical systems, and statistical mechanics. 
In this paper we formulate a general {\it geometric framework for defining physical dynamics on directed open graphs}\footnote{Note that this does \textit{not} include the (random) evolution of the graphs themselves, as studied in
random graph theory and statistical mechanics.}. The generalized Hamiltonian nature of the resulting dynamical models is due to the assumption that the constitutive relations between the variables corresponding to storage at the vertices and/or edges are derivable from an {\it energy} (Hamiltonian) function, while the remaining variables are related by static {\it energy-dissipating} relations. This will imply that the total energy itself satisfies a conservation law: the increase of the total energy is equal to externally supplied power (through the boundary vertices of the graph), minus the power lost in the dissipative elements (associated to some of the edges or vertices of the graph). The resulting generalized Hamiltonian systems, allowing for energy-dissipation and interaction with the environment, fall within the class of {\it port-Hamiltonian systems}, as coined and explored in e.g. \cite{schaftAEU95, DalsmoSIAM99, schaftspringer00, schaftGeomPhys02, Geoplex09-1}.

From a geometric point of view the generalized Hamiltonian structure of the network dynamics is defined, apart from its Hamiltonian function and energy-dissipating relations, by a {\it Dirac structure}. This Dirac structure (generalizing the symplectic or Poisson structure from classical mechanics) is directly defined by the {\it incidence} matrix of the directed graph, and thus captures the conservation laws. In fact, we will show how a directed graph gives rise to three canonically defined Dirac structures on its vertex and edge spaces. The first two of them only differ in the different role of the boundary vertices, while the third, the Kirchhoff-Dirac structure, captures the special case where no storage or dissipation is associated with the vertices of the graph (corresponding to Kirchhoff's current laws). 

We will illustrate this framework on some of the physical examples mentioned
above. Furthermore, we will show how the same port-Hamiltonian structure is shared by network dynamics from a different origin, such as consensus and clustering algorithms, and how the identification of the underlying port-Hamiltonian structure provides powerful tools for analysis and control, which unify and go beyond existing approaches.

While all examples given in the paper are simple, and could be approached from other angles as well, we believe that a major contribution of the paper resides in the identification of a common mathematical structure in all these examples, which is moreover closely related to classical Hamiltonian systems. Furthermore, the approach formalizes network dynamics as {\it open} system, and due to the compositionality properties of port-Hamiltonian systems, is easily scalable and extends to heterogeneous and multi-scale systems as well.

In a twin paper we will describe how the geometric framework
as developed in this paper for graphs can be extended to arbitrary $k$-complexes. Among others, this will allow for a structure-preserving spatial discretization of distributed-parameter physical systems, otherwise described by partial differential equations; see already \cite{schaftCDC08, schaftBosgraboek}.

Preliminary work regarding Sections 3.4 and 3.5 can be found in \cite{schaftNECSYS10, schaftCDC08, schaftBosgraboek}.

\section{From directed graphs to Dirac structures \label{sec:From-graphs-to-DS}}
\begin{figure}[htbp]
\begin{center}
\resizebox{8cm}{3cm}{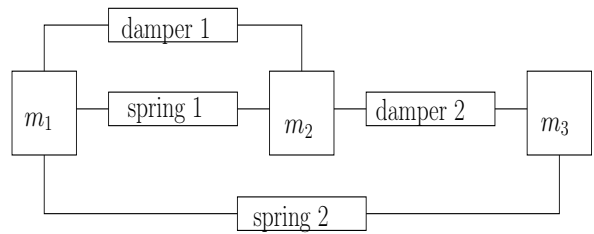}
\caption{Mass-spring-damper
      system}
\label{figure:example-network}
\end{center}
\end{figure}

\begin{figure}[htbp]
\begin{center}
\resizebox{7cm}{1.5cm}{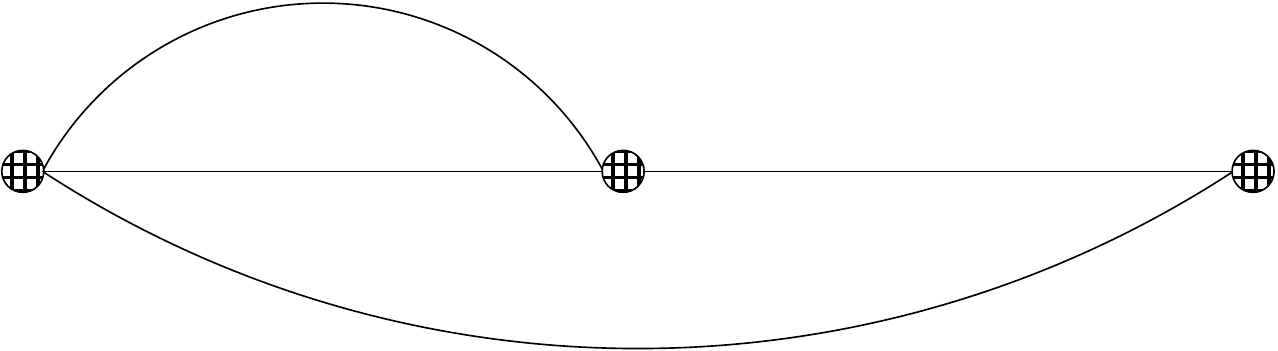}
\caption{The
      corresponding graph}
\label{figure:example-graph}
\end{center}
\end{figure}

As a guiding example let us consider a {\it mass-spring-damper} system; for example the one depicted in Figure \ref{figure:example-network}.
The underlying directed graph of such a system is defined by {\it vertices} corresponding to the masses, and {\it edges} corresponding to the springs and dampers; leading to the graph in Figure \ref{figure:example-graph}.

\noindent
How do we formalize such a system as a port-Hamiltonian system ? Key ingredient in the definition of a port-Hamiltonian system is the geometric notion of a {\it Dirac structure}, generalizing the symplectic structure from classical Hamiltonian dynamics.
In this section we will define two canonical Dirac structures on the
combination of the vertex, edge and boundary spaces of a directed graph, and their dual spaces. These two Dirac structures will only differ in the role of the boundary vertices, which for a mass-spring-damper system will be either associated to boundary masses (with inputs being the external forces on them), or will be {\it massless} (with inputs being their velocities). 

We first recall some basic notions of graph theory, see e.g. \cite{Bollobas98}, and Dirac structures, see e.g. \cite{Courant90, Dorfman93, DalsmoSIAM99}.

\subsection{Directed graphs and their vertex and edge spaces \label{sub:Directed-graphs}}

A \textit{directed graph} $\mathcal{G}= (\mathcal{V}, \mathcal{E})$ consists of a finite set $\mathcal{V}$ of \textit{vertices} (nodes) and a finite set $\mathcal{E}$ of directed \textit{edges} (branches or links), together
with a mapping from $\mathcal{E}$ to the set of ordered pairs of
$\mathcal{V}$, where no self-loops are allowed. Thus to any branch
$e\in\mathcal{E}$ there corresponds an ordered pair $(v,w)\in\mathcal{V}\times\mathcal{V}$
(with $v\not=w$), representing the tail vertex $v$ and the head
vertex $w$ of this edge.

A directed graph is completely specified by its \textit{incidence
matrix} $\hat{B}$, which is an $N\times M$ matrix, $N$ being the
number of vertices and $M$ being the number of edges, with $(i,j)$-th
element equal to $-1$ if the $j$-th edge is an edge towards vertex
$i$, equal to $1$ if the $j$-th edge is an edge originating from
vertex $i$, and $0$ otherwise. Since we will only consider directed
graphs in the sequel 'graph' will throughout mean 'directed graph'.

Given a graph, we define its \textit{vertex space} $\Lambda_{0}$
as the vector space of all functions from $\mathcal{V}$ to some linear
space $\mathcal{R}$. In the examples, $\mathcal{R}$ will be mostly
$\mathcal{R}=\mathbb{R}$ or $\mathcal{R}=\mathbb{R}^{3}$. In the
first case, $\Lambda_{0}$ can be identified with $\mathbb{R}^{N}$.
Furthermore, we define its \textit{edge space} $\Lambda_{1}$ as the
vector space of all functions from $\mathcal{E}$ to the same%
\footnote{In principle we could also associate with the edges a linear space
$\mathcal{R}^{\prime}$ which is \textit{different} from the space
$\mathcal{R}$ associated with the vertices. In that case the definition
of the incidence operator needs an additional linear map from $\mathcal{R}^{\prime}$
to $\mathcal{R}$.
} linear space $\mathcal{R}$. Again, if $\mathcal{R}=\mathbb{R}$
then $\Lambda_{1}$ can be identified with $\mathbb{R}^{M}$.

The dual spaces of $\Lambda_{0}$ and $\Lambda_{1}$ will be denoted
by $\Lambda^{0}$, respectively $\Lambda^{1}$. The duality pairing between
$f \in \Lambda_0$ and $e \in \Lambda^0$
is given as 
\[
<f\mid e>=\sum_{v\in\mathcal{V}}<f(v)\mid e(v)>\,,
\]
where $<\dot{\mid}\dot{>}$ on the right-hand side denotes the duality
pairing between $\mathcal{R}$ and $\mathcal{R}^{*}$, and a similar expression holds for $f \in \Lambda_1$ and $e \in \Lambda^1$ (with summation over the edges). 

\smallskip

The incidence matrix $\hat{B}$ of the graph induces a
linear map $B$ from the edge space to the vertex space as follows.
Define $B:\Lambda_{1}\rightarrow\Lambda_{0}$ as the linear map with
matrix representation 
$
\hat{B}\otimes I
$,
where $I:\mathcal{R}\to\mathcal{R}$ is the \textit{identity
map and} $\otimes$ denotes the Kronecker product. $B$ will be called
the \textit{incidence operator}. For $\mathcal{R}=\mathbb{R}$
the incidence operator reduces to the linear map given by the
matrix $\hat{B}$ itself, in which case we will throughout use $B$ {\it both for the incidence matrix and for the incidence operator.}
The adjoint map of $B$ is denoted as 
\[
B^{*}:\Lambda^{0}\rightarrow\Lambda^{1},
\]
and is called the \textit{co-incidence} operator. For $\mathcal{R}=\mathbb{R}^{3}$
the co-incidence operator is given by $\hat{B}^{T}\otimes I_{3}$,
while for $\mathcal{R}=\mathbb{R}$ the co-incidence operator is simply
given by the transposed matrix $\hat{B}^{T}$, and we will throughout use $B^{T}$ {\it both for the co-incidence matrix
and the for co-incidence operator.}

We will use the terminology\footnote{This terminology stems from port-based and bond-graph modeling \cite{Paynter61}, where it has a slightly more specific connotation then in our case. The space $\Lambda_0$ is also called the space of $0$-chains, while the elements of $\Lambda_1$ are called the $1$-chains. Furthermore, the dual spaces $\Lambda^0$ and $\Lambda^1$ are called the space of $0$-cochains, respectively $1$-cochains. In \cite{twin} this will be generalized to higher-order chains and co-chains. In (generalized) circuit theory, $f_1 \in \Lambda_1$ are referred to as {\it through} variables, and $fe^1 \in \Lambda^1$ as {\it across} variables.} \textit{flows} for the elements of $\Lambda_{0}$ and $\Lambda_{1}$ (notation $f_{0}$ and $f_{1}$), and \textit{efforts}
for the elements of their dual spaces $\Lambda^{0}$ and $\Lambda^{1}$
(notation $e^{0}$, respectively $e^{1}$).

\subsection{Open graphs \label{sub:Open-graphs}}

An \textit{open graph} $\mathcal{G}$ is obtained from an ordinary
graph with set of vertices $\mathcal{V}$ by identifying a subset
$\mathcal{V}_{b}\subset\mathcal{V}$ of $N_{b}$ \textit{boundary
vertices}. The interpretation of $\mathcal{V}_{b}$ is that these
are the vertices that are open to interconnection (i.e., with other
open graphs). The remaining subset $\mathcal{V}_{i}:=\mathcal{V}-\mathcal{V}_{b}$
are the $N_{i}$ \textit{internal vertices} of the open graph.

The splitting of the vertices into internal and boundary vertices
induces a splitting of the vertex space and its dual, given as
\[
\begin{array}{rcl}
\Lambda_{0} & = & \Lambda_{0i}\oplus\Lambda_{0b}\\[2mm]
\Lambda^{0} & = & \Lambda^{0i}\oplus\Lambda^{0b}\end{array}
\]
where $\Lambda_{0i}$ is the vertex space corresponding to the internal
vertices and $\Lambda_{0b}$ the vertex space corresponding to the
boundary vertices. Consequently, the incidence operator $B:\Lambda_{1}\to\Lambda_{0}$
splits as
\[
B=B_{i}\oplus B_{b}
\]
with $B_{i}:\Lambda_{1}\to\Lambda_{0i}$ and $B_{b}:\Lambda_{1}\to\Lambda_{0b}$.

Furthermore, we will define the \textit{boundary space} $\Lambda_{b}$
as the linear space of all functions from the set of boundary vertices
$\mathcal{V}_{b}$ to the linear space $\mathcal{R}$. Note that the boundary space $\Lambda_{b}$
is \textit{isomorphic} to the linear space $\Lambda_{0b}$, and that
using this isomorphism the linear mapping $B_{b}$ can be also regarded
as a mapping \[
B_{b}:\Lambda_{1}\to\Lambda_{b}\]
called the \textit{boundary incidence operator}. Nevertheless, we
will be careful in distinguishing the two isomorphic linear spaces
$\Lambda_{b}$ and $\Lambda_{0b}$ because of their different interpretations
in physical examples (e.g., for mass-spring-damper systems $\Lambda_{b}$ will denote the space of {\it external forces} as exerted on the boundary masses, and $\Lambda_{0b}$ the space of {\it momenta} of the boundary masses). The dual space of $\Lambda_{b}$ will be denoted
as $\Lambda^{b}$. The elements $f_{b}\in\Lambda_{b}$ are called
the \textit{boundary flows} and the elements $e^{b}\in\Lambda^{b}$
the \textit{boundary efforts}.

\subsection{Dirac structures \label{sub:Dirac-Structures}}

Recall (\cite{schaftAEU95,Courant90,schaftspringer00})
the definition of a (constant\footnote{This definition can be extended \cite{Dorfman93,Courant90} to (non-constant)
Dirac structures on manifolds: a Dirac structure $\mathcal{D}$ on a manifold $\mathcal{M}$ is defined as a vector subbundle of the Whitney sum $T\mathcal{M} \oplus T^*\mathcal{M}$ such that for each $m \in \mathcal{M}$ the linear space $D(m) \subset T_m\mathcal{M} \times T_m^*\mathcal{M}$ is a constant Dirac structure. This will be needed in the treatment of spatial mechanisms in Section 3.3.}) Dirac structure. Consider a vector space $\mathcal{F}$ with dual
space $\mathcal{F}^{*}$. As before, the variables $f\in\mathcal{F}$ are called
the \textit{flow} variables, while the conjugate variables $e\in\mathcal{F}^{*}$
are called the \textit{effort} variables. Define on the total space
$\mathcal{F}\times\mathcal{F}^{*}$ the \textit{indefinite inner product}
$\ll\cdot,\cdot\gg$ as 
\[
\ll(f_{a},e_{a}),(f_{b},e_{b})\gg:=<e_{a}\mid f_{b}>+<e_{b}\mid f_{a}>,f_{a},f_{b}\in\mathcal{F},e_{a},e_{b}\in\mathcal{F}^{*}
\]
where $<\cdot\mid\cdot>$ denotes the duality product between $\mathcal{F}$
and $\mathcal{F}^{*}$. \begin{definition} A subspace $\mathcal{D}\subset\mathcal{F}\times\mathcal{F}^{*}$
is a Dirac structure if 
$
\mathcal{D}=\mathcal{D}^{\pperp}\label{perp}
$, where ${}^{\pperp}$ denotes the orthogonal complement with respect
to $\ll\cdot,\cdot\gg$. \end{definition} 
In the finite-dimensional
case an equivalent, and often easier, characterization of Dirac structures
is given as follows (see e.g. \cite{Cervera07, Geoplex09-1} for a proof). 
\begin{proposition}
\label{equivalent} A subspace 
$
\mathcal{D}\subset\mathcal{F}\times\mathcal{F}^{*}
$
is a Dirac structure if and only if the following two conditions
are satisfied: 
\begin{equation}
\begin{array}{l}
(i)<e\mid f>=0,\quad \mbox{ for all }(f,e)\in\mathcal{D}\\[2mm]
(ii)\dim\mathcal{D}=\dim\mathcal{F}\end{array}\label{dirac2}
\end{equation}
\end{proposition}
Note that the first equation in (\ref{dirac2}) can be regarded as
a \textit{power-conservation property}. The second equation states that a Dirac structure has {\it maximal} dimension with respect to this power-conserving property \cite{DalsmoSIAM99, schaftspringer00}.

While Dirac structures thus formalize power-conserving interconnections of maximal dimension, the following special type of Dirac structure can be seen to be a generalization of Tellegen's theorem in circuit theory (stating that the product $<V_a \mid I_b>=0$ for {\it any} two vectors of voltages $V_a$ and currents $I_b$ satisfying Kirchhoff's laws): 
\begin{definition}
A Dirac structure $\mathcal{D}\subset\mathcal{F}\times\mathcal{F}^{*}$
is \textit{separable} if 
\begin{equation}\label{separable}
<e_{a}\mid f_{b}>=0\,,\quad \mbox{ for all } (f_{a},e_{a}),(f_{b},e_{b})\in\mathcal{D}
\end{equation}
\end{definition} 
Separable Dirac structures have the following simple geometric characterization, reminding of Kirchhoff's current and voltage laws.
\begin{proposition} Consider a separable Dirac structure $\mathcal{D}\subset\mathcal{F}\times\mathcal{F}^{*}$.
Then 
\begin{equation}\label{separable1}
\mathcal{D}=\mathcal{K}\times\mathcal{K}^{\perp}
\end{equation}
for some subspace $\mathcal{K}\subset\mathcal{F}$, where $\mathcal{K}^{\perp}= \{ e \in \mathcal{F}^* \mid <e \mid f> = 0, \forall f \in \mathcal{K} \}$. Conversely, any
subspace $\mathcal{D}$ as in (\ref{separable1}) for some subspace
$\mathcal{K}\subset\mathcal{F}$ is a separable Dirac structure. 
\end{proposition}
\textbf{Proof.} It is immediately seen that any subspace $\mathcal{K}\times\mathcal{K}^{\perp}$
satisfies (\ref{separable}), and is a Dirac structure since it satisfies
(\ref{dirac2}).
Conversely, let the Dirac structure $\mathcal{D}$ satisfy
(\ref{separable}). Define the following subspaces 
\[
\begin{array}{rcl}
\mathcal{F}_{0}=\{f\in\mathcal{F}\mid(f,0)\in\mathcal{D}\} & \, & \quad\mathcal{F}_{1}=\{f\in\mathcal{F}\mid\exists e\in\mathcal{F}^{*}\mbox{ s.t. }(f,e)\in\mathcal{D}\}\\
\mathcal{E}_{0}=\{e\in\mathcal{F}^{*}\mid(0,e)\in\mathcal{D}\} & \, & \quad\mathcal{E}_{1}=\{e\in\mathcal{F}^{*}\mid\exists f\in\mathcal{F}\mbox{ s.t. }(f,e)\in\mathcal{D}\}\end{array}
\]
It is readily seen \cite{DalsmoSIAM99} that for any Dirac structure $\mathcal{E}_{1}=(\mathcal{F}_{0})^{\perp},\mathcal{E}_{0}=(\mathcal{F}_{1})^{\perp}$.
We will now show that (\ref{separable}) implies that $\mathcal{F}_{0}=\mathcal{F}_{1}=:\mathcal{K}$
(and hence $\mathcal{E}_{0}=\mathcal{E}_{1}=:\mathcal{K}^{\perp}$).
Clearly, $\mathcal{F}_{0}\subset\mathcal{F}_{1}$. Let now $(f_{a},e_{a})\in\mathcal{D}$
and thus $f_{a}\in\mathcal{F}_{1}$. Then for all $(f_{b},e_{b})\in\mathcal{D}$
\[
\ll(f_{a},0),(f_{b},e_{b})\gg:=<e_{b}\mid f_{a}>+<0\mid f_{b}>=<e_{b}\mid f_{a}>=0
\]
by (\ref{separable}). Hence, also $(f_{a},0)\in\mathcal{D}$ and
thus $f_{a}\in\mathcal{F}_{0}$. By definition $\mathcal{F}_{0}\times\mathcal{E}_{0}\subset\mathcal{D}$,
and hence $\mathcal{K}\times\mathcal{K}^{\perp}\subset\mathcal{D}$. Finally, since the dimension of $\mathcal{K}\times\mathcal{K}^{\perp}$
equals the dimension of $\mathcal{F}$ equality results. \qed 

%
A typical instance of a separable Dirac structure, which will be frequently
used in the remainder, is the following. 
\begin{proposition}\label{prop:general}
Let $A:\mathcal{V}\to\mathcal{W}$ be a linear map between the linear
spaces $\mathcal{V}$ and $\mathcal{W}$ with adjoint mapping $A^{*}:\mathcal{W}^{*}\to\mathcal{V}^{*}$,
that is 
\begin{equation}
<w^{*}\mid Av>=<A^{*}w^{*}\mid v>\label{adjoint}
\end{equation}
for all $v\in\mathcal{V},w^{*}\in\mathcal{W}^{*}$ 
(where, as before,
$<\cdot\mid\cdot>$ denotes the duality product between the dual spaces
$\mathcal{W}$ and $\mathcal{W}^{*}$, respectively $\mathcal{V}$
and $\mathcal{V}^{*}$). Identify $(\mathcal{V} \times \mathcal{W})^* = \mathcal{V}^* \times \mathcal{W}^*$. Then 
\begin{equation}
\begin{array}{rcl}
\mathcal{D} & := & \{(v,w,v^{*},w^{*})\in(\mathcal{V}\times\mathcal{W})\times(\mathcal{V}^{*}\times\mathcal{W}^{*})\mid\\[2mm]
&  & Av=w,v^{*}= - A^{*}w^{*}\}\end{array}
\end{equation}
is a separable Dirac structure. 
\end{proposition} 
\textbf{Proof.}
Define $\mathcal{K}:=\{(v,w)\in\mathcal{V}\times\mathcal{W}\mid Av=w\}$.
Then $\mathcal{K}^{\perp}=\{(v^{*},w^{*})\in\mathcal{V}^{*}\times\mathcal{W}^{*}\mid v^{*}= - A^{*}w^{*}\}$.
\qed 
\medskip
A key feature of Dirac structures is that their {\it composition} is again a Dirac structure (in contrast with symplectic or Poisson structures, where this is not generally the case).
Let $\mathcal{D}_{A}\subset\mathcal{F}_{A}\times\mathcal{F}_{c}\times\mathcal{F}_{A}^{*}\times\mathcal{F}_{c}^{*}$
and $\mathcal{D}_{B}\subset\mathcal{F}_{B}\times\mathcal{F}_{c}\times\mathcal{F}_{B}^{*}\times\mathcal{F}_{c}^{*}$
be two Dirac structures with shared space of flow and effort variables
$\mathcal{F}_{c}$, respectively $\mathcal{F}_{c}^{*}$. Define their
\textit{composition} as 
\begin{equation}
\begin{array}{c}
\mathcal{D}_{A}\circ\mathcal{D}_{B}=\{(f_{A},e_{A},f_{B},e_{B})\in\mathcal{F}_{A}\times\mathcal{F}_{B}\times\mathcal{F}_{A}^{*}\times\mathcal{F}_{B}^{*}\mid\exists(f,e)\in\mathcal{F}_{c}\times\mathcal{F}_{c}^{*}\mbox{ s.t. }\\[2mm]
(f_{A},e_{A},f,e)\in\mathcal{D}_{A},(f_{B},e_{B},-f,e)\in\mathcal{D}_{B}\}\end{array}\label{composition}
\end{equation}
It has been shown in \cite{Cervera07,schaft99} that $\mathcal{D}_{A}\circ\mathcal{D}_{B}$ is again a Dirac structure. Separable Dirac structures turn out to have the following special compositional property: 
\begin{proposition}\label{compositional}
Let $\mathcal{D}_{A}\subset\mathcal{F}_{A}\times\mathcal{F}_{c}\times\mathcal{F}_{A}^{*}\times\mathcal{F}_{c}^{*}$
and $\mathcal{D}_{B}\subset\mathcal{F}_{B}\times\mathcal{F}_{c}\times\mathcal{F}_{B}^{*}\times\mathcal{F}_{c}^{*}$
be two separable Dirac structures given as 
\[
\mathcal{D}_{i}=\mathcal{K}_{i}\times\mathcal{K}_{i}^{\perp}\,,i=A,B,
\]
where $\mathcal{K}_{i}\subset\mathcal{F}_{i}\times\mathcal{F}_{c},i=A,B$.
Define the composition 
\[
\mathcal{K}_{A}\circ\mathcal{K}_{B}=\{(f_{A},f_{B})\in\mathcal{F}_{A}\times\mathcal{F}_{B}\mid\exists f\in\mathcal{F}_{c}\mbox{ s.t. }(f_{A},f)\in\mathcal{K}_{A},(f_{B},-f)\in\mathcal{K}_{B}\}
\]
Then the composition $\mathcal{D}_{A}\circ\mathcal{D}_{B}$ is the
separable Dirac structure 
\begin{equation}
\mathcal{D}_{A}\circ\mathcal{D}_{B}=(\mathcal{K}_{A}\circ\mathcal{K}_{B})\times(\mathcal{K}_{A}\circ\mathcal{K}_{B})^{\perp}
\end{equation}
\end{proposition}
For explicit equational representations of compositions of Dirac structures we refer to \cite{Cervera07}.

The compositionality property of Dirac structures is a key ingredient of port-Hamiltonian systems theory, and implies that the standard interconnection of port-Hamiltonian systems results in another port-Hamiltonian system with Dirac structure being the {\it composition} of the Dirac structures of the component port-Hamiltonian systems, and Hamiltonian equal to the {\it sum} of the Hamiltonians of the component systems \cite{schaft99, Cervera07}.

\subsection{The graph Dirac structures  \label{sub:Canonical-Dirac-Structures}}

We now have all ingredients to define Dirac structures corresponding to the incidence structure of a directed graph. 
\begin{definition} \label{def:GraphDiracStruc}
Consider an open graph $\mathcal{G}$ with vertex, edge
and boundary spaces, incidence operator $B$ and boundary incidence operator $B_{b}$.
The {\it flow-continuous}\footnote{The terminology {\it flow-continuous} and {\it effort-continuous} stems from the fact that in the first case the boundary flows $f_b$  are exclusively linked to the edge flows $f_1$, while in the second case the boundary efforts $e^b$ are determined by part of the internal vertex efforts $e^0$. Note that the space of involved flow and effort variables for $\mathcal{D}_{f}(\mathcal{G})$ and $\mathcal{D}_{e}(\mathcal{G})$ is {\it different}.} {\it graph Dirac structure} $\mathcal{D}_f(\mathcal{G})$ is defined as
\begin{equation}\label{flowcontinuous}
\begin{array}{rcl}
\mathcal{D}_f(\mathcal{G}) & := & \{(f_{1},e^{1},f_{0i},e^{0i},f_{b},e^{b})\in\\[2mm]
&  & \Lambda_{1}\times\Lambda^{1}\times\Lambda_{0i}\times\Lambda^{0i}\times\Lambda_{b}\times\Lambda^{b}\mid\\[2mm]
&  & B_{i}f_{1}=f_{0i},B_{b}f_{1}=f_{b},e^{1}=-B_{i}^{*}e^{0i}-B_{b}^{*}e^{b}\}\end{array}
\end{equation}
The {\it effort-continuous graph Dirac structure} $\mathcal{D}_e(\mathcal{G})$ is defined as
\begin{equation}\label{effortcontinuous}
\begin{array}{rcl}
\mathcal{D}_{e}(\mathcal{G}) & := & \{(f_{1},e^{1},f_{0},e^{0},f_{b},e^{b})\in\\[2mm]
&  & \Lambda_{1}\times\Lambda^{1}\times\Lambda_{0}\times\Lambda^{0}\times\Lambda_{b}\times\Lambda^{b}\mid\\[2mm]
&  & B_{i}f_{1}=f_{0i},B_{b}f_{1}=f_{0b}+f_{b},e^{1}=-B^{*}e^{0},e^{b}=e^{0b}\} 
\end{array}
\end{equation}
\end{definition}
By Proposition \ref{prop:general} both $\mathcal{D}_{f}(\mathcal{G})$ and $\mathcal{D}_{e}(\mathcal{G})$ are separable Dirac structures. Note that $\mathcal{D}_{f}(\mathcal{G})$ and $\mathcal{D}_{e}(\mathcal{G})$ only differ in the role of the boundary flows and efforts, and that $\mathcal{D}_{f}(\mathcal{G})= \mathcal{D}_{e}(\mathcal{G})$ if there are no boundary vertices. For mass-spring-damper systems the flow-continuous Dirac structure will correspond to the case that the boundary vertices are massless, while the effort-continuous Dirac structure corresponds to boundary masses, with momenta in $\Lambda_{0b}$.

\subsection{Interconnection of open graphs and composition of graph Dirac structures}
Interconnection of two open graphs $\mathcal{G}^{\alpha}$ and $\mathcal{G}^{\beta}$ is performed by identifying some of their boundary vertices,
and equating (up to a minus sign) the boundary efforts and flows corresponding
to these boundary vertices, resulting in a new graph. For simplicity
of exposition consider the case that the open graphs have \textit{all}
their boundary vertices in common, resulting in a (closed) graph with
set of vertices $\mathcal{V}_{i}^{\alpha}\cup\mathcal{V}_{i}^{\beta}\cup\mathcal{V}$,
where $\mathcal{V}:=\mathcal{V}_{b}^{\alpha}=\mathcal{V}_{b}^{\beta}$
denotes the set of boundary vertices of both graphs.

The incidence operator of the interconnected (closed) graph is obtained
as follows. For simplicity of notation consider the case that $\mathcal{R}=\mathbb{R}$.
Let $\mathcal{G}^{j}$ have incidence operators 
\[
B^{j}=\begin{bmatrix}B_{i}^{j}\\
B_{b}^{j}\end{bmatrix}, \quad j=\alpha,\beta 
\]
The incidence operator $B$ of the interconnected graph is then
given as 
\begin{equation}
B=\begin{bmatrix}B_{i}^{\alpha} & 0\\
0 & B_{i}^{\beta}\\
B_{b}^{\alpha} & B_{b}^{\beta}\end{bmatrix},
\end{equation}
corresponding to the interconnection constraints on the boundary
potentials and currents given by \begin{equation}
e^{b\alpha}=e^{b\beta},\quad f_{b}^{\alpha}+f_{b}^{\beta}=0\label{interc}\end{equation}
Of course, several extensions are possible. For example, one may
\textit{retain} the set of shared boundary vertices $\mathcal{V}_{b}:=\mathcal{V}_{b}^{\alpha}=\mathcal{V}_{b}^{\beta}$
as being boundary vertices (instead of internal vertices as above)
by extending (\ref{interc}) to 
\begin{equation}
e^{b\alpha}=e^{b\beta}=e^{b},\quad f_{b}^{\alpha}+f_{b}^{\beta}+f_{b}=0,\label{interc1}\end{equation}
with $f_{b},e^{b}$ the boundary flows and efforts of the interconnected
graph.

Comparing the {\it interconnection} of open graphs with the {\it composition} of their graph Dirac structures (see e.g. Proposition \ref{compositional}) it is readily seen that the flow/effort-continuous graph Dirac structure of an
interconnected graph equals the composition of the flow/effort-continuous graph Dirac structures of $\mathcal{G}^{\alpha}$ and $\mathcal{G}^{\beta}$; we leave the straightforward proof to the reader.

\subsection{Derived graph Dirac structures}\label{sub:derived}
Other Dirac structures can be {\it derived} from the flow/effort-continuous Dirac structure by {\it constraining} some of the flows and the efforts to zero. For example, the composition of the flow/effort-continuous Dirac structure with the trivial separable Dirac structure
\[
\{(f_{0i},e^{0i})\in\Lambda_{0i}\times\Lambda^{0i}\mid f_{0i}=0\}
\]
will result by Proposition \ref{compositional} in another separable Dirac structure, called the {\it Kirchhoff-Dirac structure}, which will be discussed in detail in Section \ref{sec:Kirchhoff}. 

However, there are other possibilities which we will only indicate. One, somewhat dual to the Kirchhoff-Dirac structure, is to constrain (some of) the {\it edge efforts} in the flow/effort-continuous graph Dirac structure to zero. 
Another interesting option is to constrain some of the {\it edge flows} in the flow/effort-continuous graph Dirac structure to zero. Considering the description of the flow/effort-continuous graph Dirac structure this effectively reduces (by disregarding the associated edge efforts) to the flow/effort-continuous graph Dirac structure of the reduced graph where the edges corresponding to the zero edge flows have been left out. Alternatively, one may constrain some of the {\it internal vertex efforts} to zero. Again considering the description of the flow/effort-continuous graph Dirac structure this amounts to {\it deleting} the corresponding internal vertices, turning them into boundary vertices with prescribed zero efforts. Note that this yields a setting for dealing with {\it dynamic graphs}. 

\section{Port-Hamiltonian systems on graphs \label{sec:Port-Hamiltonian-systems-on-Graphs}}

First (Section 3.1) we will describe how port-Hamiltonian systems can be defined with respect to the canonical graph Dirac structures defined above. In the subsequent subsections this will be illustrated on a number of typical examples, ranging from mass-spring-damper systems and spatial mechanisms to consensus and clustering algorithms.

\subsection{Definition of port-Hamiltonian systems with regard to the graph Dirac
structures}

In this subsection we will apply the general definition of port-Hamiltonian systems with regard to an arbitrary Dirac structure, see e.g. \cite{schaftAEU95, DalsmoSIAM99, schaftspringer00}, to the graph Dirac structures as defined above.

For clarity of exposition we throughout consider the {\it effort-continuous} graph Dirac structure $\mathcal{D}_e(\mathcal{G})$ involving the flow and effort variables
\[
(f_{1},e^{1},f_{0},e^{0},f_{b},e^{b}) \in \Lambda_{1}\times\Lambda^{1}\times\Lambda_{0}\times\Lambda^{0}\times\Lambda_{b}\times\Lambda^{b}
\]
(the exposition is directly repeated for the flow-continuous graph Dirac structure $\mathcal{D}_f(\mathcal{G})$).
A port-Hamiltonian system is specified by defining between all the {\it internal} conjugate flow and effort variables $(f_{1},e^{1},f_{0},e^{0})$ either an {\it energy-storing} relation, or a purely {\it dissipative} relation. An energy-storing relation between a vector of flow variables $f$ and a conjugate vector of effort variables $e$ is of the form\footnote{Throughout this paper $\frac{\partial H}{\partial x}(x)$ will denote the {\it column} vector of partial derivatives of $H$, with $\frac{\partial^T H}{\partial x}(x)$ denoting the {\it row} vector of partial derivatives.}
\[
\dot{x} = - f, \quad e= \frac{\partial H}{\partial x}(x),
\]
or dually
\[
\dot{x} = e, \quad f = - \frac{\partial H}{\partial x}(x),
\]
where $x$ is a vector of energy variables (of the same dimension as $f$ and $e$), and $H(x)$ is any function, representing the energy stored in the system. 

Furthermore, a dissipative relation between a vector of flow variables $f$ and a conjugate vector of effort variables $e$ is any static relation
\[
R(f,e) =0
\]
satisfying $<e \mid -f> \geq 0$ for all $(f,e)$ satisfying $R(f,e) =0$. 

In the case of a mass-spring-damper system with boundary masses (see Section \ref{sub:Mass-spring-systems}) the vertex flow and effort variables $f_{0},e^{0}$ will be related by energy-storing relations $\dot{p} = f_0, e_0=\frac{\partial K}{\partial p}(p)$, with $p$ the momenta of the masses and $K(p)$ their kinetic energies, the flow and effort variables $f_{1s},e^{1s}$ of the spring edges will correspond to energy-storing relations $\dot{q}= e^{1s}, f_{1s} = - \frac{\partial V}{\partial q}(q)$, with $q$ the spring elongations and $V(q)$ the spring potential energies, while finally the flow and effort variables $f_{1d},e^{1d}$ of the damper edges are connected by energy-dissipating relations $f_{1d} = -D(e^{1d})$ satisfying $(e^{1d})^TD(e^{1d}) \geq 0$.

\medskip

Thus a port-Hamiltonian system on a graph is defined by adding to the linear relations imposed by the graph Dirac structure {\it constitutive} relations between all the internal effort and flow variables, either of energy-storing or of dissipative type\footnote{Hence port-Hamiltonian dynamics generalizes both classical Hamiltonian dynamics (with no energy-dissipation), as well as gradient systems (where there is in general no oscillation between different energies and energy-dissipation does take place); see \cite{schaftIFAC11} and the references quoted therein.}.
It is clear that this leaves many possibilities for defining port-Hamiltonian dynamics. In particular, energy-storage, respectively dissipation, can be associated to the vertices or to the edges, or to both. The examples presented in the next subsections cover a number of these different possibilities. 

The interpretation of the flow/effort-continuous graph Dirac structure as describing {\it discrete conservation or balance laws} becomes more clear from the above description of port-Hamiltonian dynamics. For example, consider for the effort-continuous graph Dirac structure the case of energy storage associated to all the edges and vertices:
\[
\begin{array}{rclrcl}
\dot{x}^1 & = &e^1, \quad & f_1 & = &- \frac{\partial H}{\partial x^1}(x^1, x_0) \\[2mm]
\dot{x}_0 & = & -f_0, \quad & e^0 & = & \frac{\partial H}{\partial x_0}(x^1, x_0)
\end{array}
\]
for state variables $x^1 \in \Lambda^1$ and $x_0 \in \Lambda_0$, and energy function $H$.
Then the relations imposed by the effort-continuous graph Dirac structure imply
\[
\dot{x}_0 + B_if_1 =0, \quad \dot{x}^1 + B^*e^0 =0
\]
expressing discrete conservation (or balance) laws between the storage of the quantities $x_0$ associated to the vertices and the flow $f_1$ through the edges, respectively between the storage of the quantities $x^1$ associated to the edges and the effort $e^0$ at the vertices. The mass-spring system discussed in the next subsection will be of this type.

Furthermore, it is well-known \cite{schaftAEU95, DalsmoSIAM99, schaftspringer00} that port-Hamiltonian systems may easily entail {\it algebraic constraints} on their state variables. Indeed, whenever some of the effort variables $e= \frac{\partial H}{\partial x}(x)$ or $f= - \frac{\partial H}{\partial x}(x)$ are constrained by the Dirac structure, then this will generally (depending on the Hamiltonian $H$) lead to algebraic constraints on the state variables $x$.

Finally, we note a fundamental property of any port-Hamiltonian dynamics. Let $H(x)$ denote the total energy of the port-Hamiltonian system. Then because of the power-conserving property of the Dirac structure, and denoting the flows and efforts of the dissipative elements by $f_R,e^R$,
\begin{equation}
\frac{d}{dt}H(x) = <\frac{\partial^T H}{\partial x}(x) \mid \dot{x}> =  <e^R \mid f_R> + <e^b \mid f_b> \leq <e^b \mid f_b>
\end{equation}
Hence the total energy itself satisfies a conservation law: its increase is equal to the externally supplied power $<e^b \mid f_b>$ minus the dissipated power $-<e^R \mid f_R>$.

\begin{remark}
One may directly extend the definition of port-Hamiltonian systems on graphs to the case where the graphs are dynamically changing in time; as briefly indicated in Section \ref{sub:derived}. This leads to {\it switching port-Hamiltonian systems on graphs}; see already \cite{gerritsen, schaftcamlibel, Geoplex09-1}.
\end{remark}

\subsection{Mass-spring-damper systems \label{sub:Mass-spring-systems}}
The basic way of modeling a \textit{mass-spring-damper system} as a port-Hamiltonian system on a graph $\mathcal{G}$ is
to associate the \textit{masses} to the \textit{vertices}, and the
\textit{springs} and \textit{dampers} to the \textit{edges} of the graph, cf. Figures \ref{figure:example-network, figure:example-graph}.
For clarity of exposition we will start with the separate treatment of mass-spring (Section \ref{sub:massspring}) and mass-damper (Section \ref{sub:massdamper}) systems, before their merging in Section \ref{sub:massspringdamper}.

\subsubsection{Mass-spring systems}\label{sub:massspring}
Consider a graph $\mathcal{G}$ with $N$ vertices (masses)
and $M$ edges (springs), specified by an incidence operator $B$.
First consider the situation that the mass-spring system is
located in one-dimensional space $\mathcal{R}=\mR$, and the springs
are scalar. A vector in the vertex space $\Lambda_{0}$ then corresponds to the
vector $p$ of the scalar momenta of all $N$ masses, i.e., $p\in\Lambda_{0}=\mR^{N}$.
Furthermore, a vector in the dual edge space $\Lambda^{1}$ will
correspond to the total vector $q$ of elongations of all $M$ springs,
i.e., $q\in\Lambda^{1}=\mR^{M}$.
Next ingredient is the definition of the Hamiltonian (stored energy)
$H: \Lambda^{1}\times\Lambda_{0}\to\mathbb{R}$ (which normally splits into a sum of the kinetic and potential energies
of each mass and spring).
In the absence of boundary vertices the dynamics of the mass-spring
system is then described as the port-Hamiltonian system 
\begin{equation}
\begin{bmatrix}\dot{q}\\
\dot{p}\end{bmatrix}=\begin{bmatrix}0 & B^{T}\\
-B & 0\end{bmatrix}\begin{bmatrix}\frac{\partial H}{\partial q}(q,p)\\
\frac{\partial H}{\partial p}(q,p)\end{bmatrix}
\end{equation}
defined with respect to the graph Dirac structure $\mathcal{D}_{e}(\mathcal{G}) = \mathcal{D}_{f}(\mathcal{G})$. Note that in fact the skew-symmetric matrix 
\begin{equation}\label{Poisson}
J:= \begin{bmatrix}0 & B^{T}\\
-B & 0\end{bmatrix}
\end{equation}
defines a \textit{Poisson structure} on the state space $\Lambda^{1}\times\Lambda_{0}$.

The inclusion of boundary vertices, and thereby of external interaction,
can be done in different ways. The first option is to associate \textit{boundary
masses} to the boundary vertices. Considering the effort-continuous graph Dirac
structure $\mathcal{D}_{e}(\mathcal{G})$ we are then led to the port-Hamiltonian
system 
\begin{equation}\label{mass-spring1}
\begin{array}{rcl}
\dot{q} & = & B^{T}\frac{\partial H}{\partial p}(q,p)\\[2mm]
\dot{p} & = & -B\frac{\partial H}{\partial q}(q,p)+ Ef_{b}\\[2mm]
e^{b} & = & E^T\frac{\partial H}{\partial p}(q,p)
\end{array}
\end{equation}
Here $E$ is a matrix with as many columns as there are boundary vertices; each column consists of zeros except for exactly one $1$ in the row corresponding to the associated boundary vertex.
$f_{b}\in\Lambda_{b}$ are the external \textit{forces}
exerted (by the environment) on the boundary masses, and $e_{b}\in\Lambda^{b}$
are the \textit{velocities} of these boundary masses.

Another possibility is to start from the flow-continuous graph Dirac structure
$\mathcal{D}_f(\mathcal{G})$. In this case there are no masses associated
to the boundary vertices, and we obtain the port-Hamiltonian system
(with $p$ now denoting the vector of momenta of the masses associated
to the \textit{internal} vertices)
\begin{equation}\label{mass-spring2}
\begin{array}{rcl}
\dot{q} & = & B_{i}^{T}\frac{\partial H}{\partial p}(q,p)+B_{b}^{T}e^{b}\\[2mm]
\dot{p} & = & -B_{i}\frac{\partial H}{\partial q}(q,p)\\[2mm]
f_{b} & = & B_{b}\frac{\partial H}{\partial q}(q,p)\end{array}
\end{equation}
with $e^{b}\in\Lambda^{b}$ the velocities of the massless boundary
vertices, and $f_{b}\in\Lambda_{b}$ the forces at the boundary vertices
as \textit{experienced} by the environment.
Note that in this latter case the external velocities
$e^{b}$ of the boundary vertices can be considered to be \textit{inputs}
to the system and the forces $f_{b}$ to be \textit{outputs}; in contrast
to the previously considered case (boundary vertices corresponding to boundary
masses), where the forces $f_{b}$ are inputs and the velocities $e^{b}$
the outputs of the system\footnote{One can also consider the hybrid case where \textit{some} of
the boundary vertices are associated to masses while the remaining
ones are massless.}. 

The above formulation of mass-spring systems in $\mathcal{R}=\mR$
directly extends to $\mathcal{R}=\mR^{3}$ by using the incidence
operator $B=\hat{B}\otimes I_{3}$ as defined before.
Finally, we remark that in the above treatment we have considered springs with {\it arbitrary} elongation vectors $q \in \Lambda^1$. For ordinary springs the vector $q$ of elongations is given as $q= B^Tq_c$, where $q_c \in \Lambda^0$ denotes the vector of positions of the vertices. Hence in this case $q \in \im B^T \subset \Lambda^1$. Note that the subspace $\im B^T \times \Lambda_0 \subset \Lambda^1 \times \Lambda_0$ is an invariant subspace with regard to the dynamics (\ref{mass-spring1}) or (\ref{mass-spring2}). We will return to this in Section \ref{subsec:symmetry}.

\subsubsection{Mass-damper systems \label{sub:massdamper}}

Replacing springs by dampers leads to \textit{mass-damper systems}. In the case of the flow-continuous graph Dirac structure this yields the following equations 
\begin{equation}\label{massdamper}
\begin{array}{rcl}
B_{i}f_{1} & = & -\dot{p}\\[2mm]
B_{b}f_{1} & = & f_{b}\\[2mm]
e^{1} & = & -B_{i}^{T}\frac{\partial H}{\partial p}(p)-B_{b}^{T}e^{b}\end{array}\end{equation}
where $f_{1},e^{1}$ are the flows and efforts corresponding to the
dampers (damping forces, respectively, velocities). For example, for
{\it linear} dampers $f_{1}=-Re^{1},$
where $R$ is the positive diagonal matrix with the damping constants on its diagonal. Substitution
into (\ref{massdamper}) then yields the port-Hamiltonian system
\begin{equation}\label{mass-damper}
\begin{array}{rcl}
\dot{p} & = & -B_{i}RB_{i}^{T}\frac{\partial H}{\partial p}(p)-B_{i}RB_{b}^{T}e^{b}\\[2mm]
f_{b} & = & B_{b}RB_{i}^{T}\frac{\partial H}{\partial p}(p)+B_{b}^{T}RB_{b}^{T}e^{b}
\end{array}
\end{equation}
where, as before, the inputs $e^{b}$ are the boundary velocities
and $f_{b}$ are the forces as experienced at the massless boundary
vertices. 

\subsubsection{Mass-spring-damper systems}\label{sub:massspringdamper}
For a mass-spring-damper system the edges will correspond partly to springs, and partly to dampers.
Thus a mass-spring-damper system is described by a graph $\mathcal{G} (\mathcal{V}, \mathcal{E}_s \cup \mathcal{E}_d)$, where the vertices in $\mathcal{V}$ correspond to the {\it masses}, the edges in $\mathcal{E}_s$ to the {\it springs}, and the edges in $\mathcal{E}_d$ to the dampers of the system. This corresponds to an incidence matrix $B= \begin{bmatrix} B_s & B_d \end{bmatrix}$, where the columns of $B_s$ reflect the spring edges and the columns of $B_d$ the damper edges. For the case {\it without} boundary vertices the dynamics of such a mass-spring-damper system with linear dampers takes the form
\begin{equation}\label{massspringdamper0}
\begin{bmatrix} \dot{q} \\ \dot{p} \end{bmatrix} = \begin{bmatrix} 0 & B_s^T \\ -B_s & -B_dRB_d^T \end{bmatrix} \begin{bmatrix} \frac{\partial H}{\partial q}(q,p) \\ \frac{\partial H}{\partial p}(q,p) \end{bmatrix}
\end{equation}
In the presence of boundary vertices we may distinguish, as above, between {\it massless} boundary vertices, with inputs being the boundary velocities and outputs the boundary (reaction) forces, and {\it boundary masses}, in which case the inputs are the external forces and the outputs the velocities of the boundary masses. We leave the details to the reader.

\subsection{Spatial mechanisms\label{sub:Spatial-mechanism}}

In this section we briefly discuss the extension of mass-spring-damper systems in $\mathbb{R}$ or $\mathbb{R}^3$ to {\it spatial mechanisms}, that is, networks of rigid bodies in $\mathbb{R}^3$ related by joints. In this
case, the linear space $\mathcal{R}$ is given by $\mathcal{R}:=\mathrm{se}^{*}(3)$,
the dual of the Lie algebra of the Lie group $SE(3)$ describing the
position of a rigid body in $\mR^{3}$. A spatial mechanism (or \emph{multibody system}) is a mechanical system
consisting of \emph{rigid bodies} related by joints (defined as \emph{kinematic pairs}) restricting the relative motion between
the rigid bodies. The reader may find numerous references about their
definition and analysis in \cite{Ball00} \cite{Selig96}, using
different geometric representations of rigid body displacements. In
this paper however we shall follow the exposition in e.g. \cite{murray94, Karger_Novak85},
which is based on the Lie group of isometries in Euclidean space
$\mathbb{R}^{3}$.

The basic topology of the mechanism is described by a directed graph, called the \emph{primary graph},
whose vertices correspond to the rigid bodies and whose edges are
associated with the kinematic pairs. This is similar to the mass-spring or mass-damper systems described in Section \ref{sub:Mass-spring-systems}, with the difference
that the dynamical system associated with each vertex corresponds
to rigid body dynamics instead of point-mass dynamics, and that the edges are in first instance associated with kinematic
constraints between the bodies instead of springs or dampers. We shall see how (spatial) springs
may be included in second instance.


\subsubsection{The rigid body element}
\label{rigidbody} 

The configuration space of a rigid body is the Lie group of isometries
in Euclidean space $\mathbb{R}^{3}$, called the Special Euclidean Group
and denoted by $SE(3)\ni Q$ (also called the space of \emph{rigid body
displacements}). Using the momentum map associated with the action
of $SE\left(3\right)$ on its cotangent bundle $T^{*}SE\left(3\right)$,
following for instance \cite[chap. 4]{Libermann_marle87}, one may
define the state space of the rigid body as $SE\left(3\right)\times se^{*}\left(3\right)\ni\left(Q,\, P\right)$
by means of the left trivialization, where $P$ is called the \emph{momentum
in body frame}.

The kinetic energy of a rigid body is defined by
\begin{equation}
K(P)=\frac{1}{2}\langle P,\,\left(I^{\flat}\right)^{-1}(P)\rangle\label{maschkeeqkin_energy}\end{equation}
where $I^{\flat}:\: se\left(3\right)\rightarrow se^{*}\left(3\right)$
is is a symmetric, positive-symmetric isomorphism, called the \emph{inertia
operator of the rigid body in the body frame.}The potential energy of the
rigid body is defined by a function $U(Q)$ of the displacement $Q$.
The potential energy may be due to gravity or may
be zero in the case of the Euler-Poinsot problem.

We assume that the rigid body is subject to an external force expressed as an element $W_{e}\in se^{*}\left(3\right)$,
called \emph{force in fixed frame} \cite{Libermann_marle87} or \emph{wrench in fixed frame} \cite{Karger_Novak85}), which is obtained
by the right trivialization of $T^{*}SE\left(3\right)$. We shall
associate a conjugate velocity to this external force, the \emph{velocity
}of the body\emph{ $T_{e}$ in fixed frame} \cite{Libermann_marle87}
(also called \emph{twist in fixed frame} \cite{Karger_Novak85}),
and obtained by the right trivialization of $TSE\left(3\right)$.

The dynamical equations of the rigid body elements may then be written
as a port-Hamiltonian system \cite{schaftAEU95} \cite[eqn. (1.37)]{maschkeLondon97_tot}:
\begin{equation}
\begin{array}{l}
\frac{d}{dt}\left(\begin{array}{c}
Q\\
P\end{array}\right)=\left(\begin{array}{cc}
0 & TL_{Q}\\
-T^{*}L_{Q} & -P\times\end{array}\right)\,\left(\begin{array}{c}
dU(Q)\\
\left(I^{\flat}\right)^{-1}(P)\end{array}\right)+\left(\begin{array}{c}
0\\
Ad_{Q}^{*}\end{array}\right)\, W_{e}\\
T_{e}=\left(\begin{array}{cc}
0 & Ad_{Q}\end{array}\right)\,\left(\begin{array}{c}
dU(Q)\\
\left(I^{\flat}\right)^{-1}(P)\end{array}\right)\end{array}\label{maschkeeqconst_body}
\end{equation}
where $TL_{Q}$ denotes the tangent map to the \emph{left translation}
(mapping the velocities $T\in se\left(3\right)$ in body frame into
the velocities $v\in T_{Q}SE\left(3\right)$), $T^{*}L_ {Q}$ denotes
its dual map (mapping forces $F\in T_{Q}^{*}SE\left(3\right)$ into
forces in body frame $W\in se^{*}\left(3\right)$), $Ad_{Q}$ denotes
the \emph{adjoint representation} (mapping velocities in body frame
into velocities in fixed frame), $Ad_{Q}^{*}$ denotes the adjoint
map to $Ad_{Q}$, while finally  $\times$ is defined by \emph{the coadjoint representation
of the Lie algebra }$se\left(3\right)$, that is, $W\times T=ad_{T}^{*}W$, for any $\left(W,\, T\right)\in se^{*}\left(3\right)\times se\left(3\right)$.
The Dirac structure $\mathcal{D}_{RB}$ of this port-Hamiltonian system (\ref{maschkeeqconst_body}) is thus specified as\footnote{Note that this is a non-constant Dirac structure on $SE(3)$.}
\begin{equation}\label{eq:DS_Rigid_Body}
\begin{array}{l}
\mathcal{D}_{RB}(Q) =  \left\{ \left(v,\, W,\, T_{e},\, F,\, T,\, W_{e}\right) \in \right. \\[2mm]
T_{Q}SE(3)\times se^{*}(3)\times se(3)\times T_{Q}^{*}SE(3)\times se(3)\times se^{*}(3)\,\mid \\[2mm]
\left(\begin{array}{c}
v\\
W\end{array}\right)=\left(\begin{array}{cc}
0 & TL_{Q}\\
-T^{*}L_{Q} & -P\times\end{array}\right)\,\left(\begin{array}{c}
F\\
T\end{array}\right)+\left(\begin{array}{c}
0\\
Ad_{Q}^{*}\end{array}\right)\, W_{e} \, ,\\
\left.T_{e}=\left(\begin{array}{cc}
0 & Ad_{Q}\end{array}\right)\,\left(\begin{array}{c}
F\\
T\end{array}\right)\right\} \end{array}
\end{equation}
In this way we have associated with every vertex of the primary graph of the spatial mechanism a dynamical system (\ref{maschkeeqconst_body}) with inputs and outputs $\left(W_{e},\, T_{e}\right)\in se^{*}\left(3\right)\times se\left(3\right)$.

\subsubsection{The kinematic pair}\label{kinematicpair}

Constraints between the rigid bodies of the mechanism will be specified by \emph{kinematic pairs} corresponding to each edge of the primary graph.
A \emph{kinematic pair} is the 
idealization of a set of contacts that occur between two rigid bodies
at some configuration of the bodies. It constrains the possible relative twists between the bodies as well as the possible transmitted wrenches.
The wrench $W$ transmitted by
a kinematic pair is constrained to a linear subspace of the space
of wrenches $se^{*}(3)$ called the \emph{space of constraint wrenches}
and denoted by ${\mathcal{CW}}$. A relative twist between the two
bodies is allowed by the kinematic pair when it produces no work with
any transmissible wrench. The relative twist is thus constrained to a linear subspace ${\mathcal{FT}}$ of the space of twists $se(3)$, called
the \emph{space of freedom twists}. Since an ideal kinematic pair is workless the subspace ${\mathcal{FT}}$ is orthogonal (in the sense of the duality product) to the space of transmitted wrenches ${\mathcal{CW}}$, that is
$\mathcal{FT} = \mathcal{CW}^{\perp}$.

We have defined the spaces of freedom twists and constraint wrenches as subspaces of of the Lie algebra $se(3)$ and its dual. However these spaces express constraints on the twists and wrenches of the rigid bodies related by the kinematic pairs and hence are expressed in some common frame with configuration $Q_{KP}$. (In most cases equal to the configuration of one of the related bodies.) Consequently, the constitutive relations of a kinematic pair are given in terms of its pair of twists and wrenches $\left(T_{KP},\,W_{KP}\right) \in T_{Q_{KP}}SE\left(3\right) \times T^*_{Q_{KP}}SE\left(3\right)$ in the form
\begin{equation}
Ad_{Q_{KP}}^{*}\:W_{KP} \in \mathcal{CW} \quad \textrm{and} \quad Ad_{Q^{-1}_{KP}}\:T_{KP} \in \mathcal{FT}
\end{equation}
Hence the constitutive equations of a kinematic pair may be expressed as the following non-constant separable Dirac structure:
\begin{equation}\label{maschkeeqdiracKP}
\begin{array}{rcl}
{\mathcal{D}}_{\mathcal{CW}}\left(Q_{KP}\right)&=&\left\{\left(T_{KP},\,W_{KP}\right) \in T_{Q_{KP}}SE\left(3\right) \times T^*_{Q_{KP}}SE\left(3\right), \mid \right. \\[2mm]
& & \left. Ad_{Q_{KP}}^{*}\:W_{KP} \in \mathcal{CW} \quad \textrm{and} \quad Ad_{Q^{-1}_{KP}}\:T_{KP} \in  \mathcal{CW} ^{\perp} \right\}
\end{array}
\end{equation} 
The kinematic pair introduced above represents ideal kinematic constraints;
in general, however, mechanical work may be produced at
the kinematic pair due to the presence of actuators or springs and dampers. Such an interaction is captured by considering the linear space ${\mathcal{IW}}:=se^{*}(3)/{\mathcal{CW}}$ (which may be identified with a subspace of $se^{*}(3)$ complementary to the space
of constraint wrenches $\mathcal{CW}$). The space of interaction twists is then defined
as its dual space ${\mathcal{IT}}:={\mathcal{IW}}^*\simeq \mathcal{CW}^{\perp}$. Using the canonical projection $\pi$
of $se^{*}(3)$ onto ${\mathcal{IW}}$, together with its adjoint
map ${\pi}^{*}$, one may thus define an additional pair of port variables enabling to connect actuators, damper or spring elements to the kinematic pairs. 
The resulting \emph{interacting} kinematic pair is then defined as a $2$-port element with constitutive relations defined by the following non-constant separable Dirac structure 
\begin{equation}
\label{maschkeeqdiracIKP}
\begin{array}{rcl}
{\mathcal{D}}^I_{\mathcal{CW}}\left(Q_{KP}\right) &= &\left\{\left(T_{KP},\,W_{KP},\,T_I,\,W_I\right)  \in \right. \\[3mm]
&&T_{Q_{KP}}SE(3) \times T^*_{Q_{KP}}SE(3) \times \mathcal{CW}^{\perp} \times se^{*}(3)/{\mathcal{CW}} \mid \\[3mm]
&& \left. W_I = \pi \circ Ad_{Q_{KP}}^{*} \left(W_{KP} \right) \, , T_{KP}=-Ad_{Q_{KP}} \circ \pi ^*\left(T_I \right) \right\}
\end{array}
\end{equation} 
It is easy to check that for $W_I=0$ the interacting kinematic pair reduces to the kinematic pair as defined before.

\subsubsection{The kinestatic connection network \label{kinestatic}}

The primary graph of the mechanism together with the kinematic pairs is called the \emph{kinestatic model} of the mechanical system. Its associated Dirac structure is the \emph{composition} of the Dirac structures corresponding to the kinematic pairs with the flow-continuous\footnote{Or the effort-continuous graph Dirac structure in case the rigid bodies corresponding to the boundary vertices have non-zero inertia operator.} graph Dirac structure of the primary graph.

Consider a mechanism defined by its primary graph ${\mathcal{G}}$
composed of $n_{RB}$ \emph{internal} vertices (associated with the rigid bodies), $n_{b}$
boundary vertices corresponding to rigid bodies with zero inertia operator
and $n_{KP}$ edges (associated with the kinematic pairs). Define
the vertex space $\Lambda_{0}\ni T^{RB},$ and the edge space
$\Lambda_{1}\ni T^{KP}$ with respect to the Lie algebra $se\left(3\right)$, which represent respectively the external twist of the rigid bodies
and the kinematic pairs. The dual spaces $\Lambda^{0}\ni W^{RB}$
, respectively $\Lambda^{1}\ni W^{KP}$, then represent the external wrenches of the rigid bodies, respectively the wrenches
of the kinematic pairs; see also Remark 2.1. The twists and wrenches of the
boundary vertices (the rigid bodies with zero inertia operator) are associated with the vertex space $\Lambda_{b}\ni T^{b}$,
respectively its dual $\Lambda^{b} \ni W^{b}$. Kirchhoff's
laws on the twists and wrenches \cite{Davies83} amount to constraining these variables to belong to the flow-continuous graph Dirac
structure, i.e.,
\[
\left( T^{KP},\, W^{KP},\, T^{RB},\, W^{RB},\, T^{b},\, W^{b}\right) \in\mathcal{D}_f(\mathcal{G})
\]
Composition of $\mathcal{D}_f(\mathcal{G})$ with the Dirac structures $\mathcal{D}_{\mathcal{CW}}(Q_{KP})$ corresponding to all the kinematic pairs then results in the Dirac structure $\mathcal{D}_{KS}$ of the kinestatic model:
%
\begin{equation} \label{eq:Kinestatic_DS}
\left( T^{I},\, W^{I},\, T^{RB},\, W^{RB},\, T^{b},\, W^{b}\right) \in { \mathcal{D}}_{KS}
\end{equation}

\subsubsection{Dynamics of spatial mechanisms}
\label{dynamics} \indent
The state space $\mathcal{X}$ of the complete mechanism is the product space of the state spaces of all the rigid bodies, i.e., $\mathcal{X}= \left(SE(3) \times se^{*}(3)\right)^{n_{RB}}$,
where $n_B$ denotes the number of rigid bodies (equal to the number of internal vertices of the primary graph).
Recalling that the rigid body dynamics is defined as a port Hamiltonian system with respect to the Dirac structure (\ref{eq:DS_Rigid_Body}) one then obtains the overall Dirac structure $\mathcal{D}_M$ of the mechanism by composing the Dirac structure ${ \mathcal{D}}_{KS}$ of the kinestatic model with the Dirac structures ${ \mathcal{D}}_{RB}$ of all the rigid bodies.
Finally, defining the Hamiltonian $H_M\left(x\right)$ as the \emph{sum} of the Hamiltonians of each body one obtains the following port-Hamiltonian model of the mechanism: 
\begin{equation} \label{eq:PHS_Mechanism}
\left(- \frac{dx}{dt},\, \frac{\partial H_M}{\partial x}(x),\, T^{I},\, W^{I},\, T^{b},\, W^{b}\right) \in \mathcal{D}_{M} 
\end{equation}

\subsection{Hydraulic networks}
The interpretation of the flow-/effort-continuous graph Dirac structures as capturing the basic conservation/balance laws of a network becomes especially tangible for hydraulic networks.

A hydraulic network can be modeled as a directed graph with edges corresponding to pipes, see e.g. \cite{bookhydraulic, DePersis}. The vertices may either correspond to connection points with {\it fluid reservoirs} (buffers), or merely to connection points of the pipes; we concentrate on the first case (the second case corresponding to a Kirchhoff-Dirac structue, cf. Section \ref{sub:The-Kirchhoff-Dirac-structure}).
Let $x_v$ be the stored fluid at vertex $v$ and let $\nu_e$ be the fluid flow through edge $e$. Collecting all stored fluids $x_v$ into one vector $x$, and all fluid flows $\nu_e$ into one vector $\nu$, the {\it mass-balance} is summarized as
\begin{equation}
\dot{x} = - B\nu
\end{equation}
with $B$ denoting the incidence matrix of the graph. In the absence of fluid reservoirs this simply reduces to Kirchhoff's current laws $B\nu=0$. 

For incompressible fluids a standard model of the fluid flow $\nu_e$ through pipe $e$ is
\begin{equation}
J_e \dot{\nu}_e = P_i - P_j - \lambda_e(\nu_e)
\end{equation}
where $P_i$ and $P_j$ are the pressures at the tail, respectively head, vertices of edge $e$. Note that this captures in fact {\it two} effects; one corresponding to energy storage and one corresponding to energy dissipation. Defining the energy variable $\varphi_e := J_e \nu_e$ the stored energy in the pipe associated with edge $e$ is given as $\frac{1}{2J_e} \varphi_e^2 = \frac{1}{2}J_e \nu_e^2$. Secondly, $\lambda_e(\nu_e)$ is a damping force corresponding to energy dissipation.

In the case of fluid reservoirs at the vertices the pressures $P_v$ at each vertex $v$ are functions of $x_v$, and thus, being scalar functions, always derivable from an energy function
$
P_v = \frac{\partial H_v}{\partial x_v}(x_v), v \in \mathcal{V},
$
for some Hamiltonian $H_v(x_v)$ (e.g. gravitational energy). The resulting dynamics (with state variables $x_v$ and $\varphi_e$) is port-Hamiltonian with respect to the graph Dirac structure $\mathcal{D}_f(\mathcal{G}) = \mathcal{D}_e(\mathcal{G})$.  The set-up is immediately extended to boundary vertices (either corresponding to controlled fluid reservoirs or direct in/ or out-flows). 

\subsection{Port-Hamiltonian formulation of consensus algorithms \label{sub:Consensus-algorithm}}

While all previous examples of port-Hamiltonian systems on graphs
arise from physical modeling, the system treated in this subsection has a different
origin. Nevertheless, it shares the same structure, and in fact, turns out to have dynamics equal to the mass-damper system treated before.

Consider a network of $N$ agents moving in linear space $\mathcal{R}$,
whose interaction topology is described by an \textit{undirected}
graph $\mathcal{G}$ (symmetric interaction). Denote by $E(\mathcal{G})$
the edges of this undirected graph, consisting of unordered pairs
$(v,w)$ of vertices $v,w$. Hence $(v,w)\in E(\mathcal{G})$ if and
only if $(w,v)\in E(\mathcal{G})$. Thus the vertices of the graph
correspond to the agents, and the edges to the symmetric interactions
between them. Distinguish between \textit{leader} and \textit{follower}
agents, see e.g. \cite{Rahmani09}, and associate the leader agents
to the boundary vertices and the follower agents to the internal vertices.

Associated to each agent $v$ there is a vector $x_{v}\in\mathcal{R}$
describing the motion in the linear space $\mathcal{R}$. In the standard
consensus algorithm, see e.g. \cite{Olfati-SaberProcIEEE07}, the vector
$x_{v}$ of each follower agent $v$ satisfies the following dynamics
\begin{equation}\label{consensus}
\dot{x}_{v}(t)=-\sum_{(v,w)\in E(\mathcal{G})}g_{(v,w)}(x_{v}(t)-x_{w}(t))
\end{equation}
where $g_{(v,w)}>0$ denotes a certain positive-definite \textit{weight}
matrix associated to each edge.
For simplicity of exposition let us take the linear space $\mathcal{R}$
to be equal to $\mathbb{R}$ in the rest of this section, implying
that $g_{(v,w)}>0$ are just positive constants. Collecting all \textit{follower}
variables $x_{v}$ into one vector $x\in\mathbb{R}^{N_{i}}$, and
all \textit{leader} variables $x_{v}$ into one vector $u\in\mathbb{R}^{N_{b}}$,
it is readily checked that the dynamics can be written as 
\begin{equation}\label{consensus1}
\dot{x}=-B_{i}GB_{i}^{T}x-B_{i}GB_{b}^{T}u
\end{equation}
with $B$ the incidence matrix of the graph \textit{endowed with
an arbitrary orientation}\footnote{It is easily seen \cite{Bollobas98} that the Laplacian matrix $BGB^T$ is {\it independent} of the chosen orientation.}, and $G$ the diagonal matrix with elements
$g_{(v,w)}$ corresponding to each edge $(v,w)$. This defines a port-Hamiltonian
system with respect to the flow-continuous graph Dirac structure $\mathcal{D}_f(\mathcal{G})$ and the Hamiltonian $H(x):=\frac{1}{2}\parallel x\parallel^{2}$. Indeed, (\ref{consensus1}) is equal to 
\begin{equation}\label{PHconsensus}
\dot{x}=-B_{i}GB_{i}^{T}\frac{\partial H}{\partial x}(x)-B_{i}GB_{b}^{T}u,
\end{equation}
which are the same equations as for the mass-damper system (\ref{mass-damper}),
with $u=e_{b}\in\Lambda^{b}$.
%
Note that the corresponding artificial \textit{output} vector $y=f_{b}\in\Lambda_{b}$ given as
\[
y := B_{b}GB_{i}^{T}\frac{\partial H}{\partial x}(x)+B_{b}GB_{b}^{T}u
\]
equals minus the rate of the leader variables if the leader variables
were supposed to obey the consensus algorithm with regard to the follower
agents (which is \textit{not} the case). Hence this artificial output
measures the discrepancy between the leaders and the followers.

\subsubsection{Network clustering dynamical models}\label{sub:clustering}

A dynamical model for network clustering, where the network splits into subnetworks which separately reach consensus, was recently proposed and discussed in \cite{burger}. Consider again a multi-agent system with $N$ agents and state variables $x_i \in \mathbb{R}, i=1,\cdots,N$, whose dynamics is described as
\begin{equation}
\dot{x}_i = - \frac{dJ_i}{dx_i}(x_i) + u_i, \quad i=1,\cdots,N
\end{equation}
where the functions $J_i(x_i)$ are certain objective functions. Let the vector $u$ with components $u_i$ be determined as
\begin{equation}
u  =  B\frac{\partial V}{\partial z}(z), \quad
\dot{z}  =  -B^Tx
\end{equation}
where $V(z)=V_1(z_1) + \cdots + V_M(z_M)$ for certain functions $V_1,\cdots,V_M$.
This is readily seen to result in a port-Hamiltonian system with total Hamiltonian $H(x,z) = \frac{1}{2} \|x\|^2 + V(z)$, and a nonlinear resistive characteristic associated to each $i$-th vertex defined by the functions $J_i(x_i)$, interpreted as Rayleigh dissipation functions\footnote{The condition of convexity imposed in \cite{burger} on $J_i$ thus corresponds to incremental passivity.}. {\it Clustering} may occur once the energy functions $V_i$ define {\it bounded} constitutive relations $e_{1i} = \frac{dV_i}{dz_i}(z_i)$ for the edge efforts. Depending on the strength of the objective functions $J_i$ this will imply that consensus among the $x_i$-variables will only be reached for subnetworks.

Many other models of network dynamics, of a 'non-physical' background, can be formulated as port-Hamiltonian systems on graphs. Examples include coordination control \cite{Arcak07} and edge agreement \cite{zelazo}.

\section{Dynamical analysis}\label{sec:analysis}

In this section we will investigate the dynamical properties of a paradigmatic example of a port-Hamiltonian system on a graph, namely the mass-spring-damper system as discussed in Section \ref{sub:massspringdamper}. As we have seen, many other examples share the same mathematical structure, and the dynamical analysis for other examples will follow the same lines.

Thus we will consider a mass-spring-damper system as described by a graph $\mathcal{G} (\mathcal{V}, \mathcal{E}_s \cup \mathcal{E}_d)$, where the vertices in $\mathcal{V}$ correspond to the {\it masses}, the edges in $\mathcal{E}_s$ to the {\it springs}, and the edges in $\mathcal{E}_d$ to the dampers of the system, with incidence matrix $B= \begin{bmatrix} B_s & B_d \end{bmatrix}$, where the columns of $B_s$ reflect the spring edges and the columns of $B_d$ the damper edges. Without boundary vertices the dynamics takes the form (see equation (\ref{massspringdamper0}) in Section \ref{sub:massspringdamper})
\begin{equation}\label{massspringdamper}
\begin{bmatrix} \dot{q} \\ \dot{p} \end{bmatrix} = \begin{bmatrix} 0 & B_s^T \\ -B_s & -B_dRB_d^T \end{bmatrix} \begin{bmatrix} \frac{\partial H}{\partial q}(q,p) \\ \frac{\partial H}{\partial p}(q,p) \end{bmatrix}
\end{equation}
Throughout this section we make the simplifying assumption\footnote{This assumption can be made without loss of generality, since otherwise the same analysis can be repeated for each connected component.}:
\begin{assumption}
The graph $\mathcal{G} (\mathcal{V}, \mathcal{E}_s \cup \mathcal{E}_d)$ is connected, or equivalently $\ker B_s^T \cap \ker B_d^T = \spa \mathds{1}$.
\end{assumption}

\subsection{Equilibria and Casimirs}
\begin{proposition}
The set of equilibria $\mathcal{E}$ of (\ref{massspringdamper}) is given as
\begin{equation}
\mathcal{E} = \{ (q,p) \in \Lambda^1 \times \Lambda_0 \mid \frac{\partial H}{\partial q}(q,p) \in \ker B_s,
\frac{\partial H}{\partial p}(q,p) \in   \spa \mathds{1}\}
\end{equation}
\end{proposition}
\noindent
{\bf Proof.}
$(q,p)$ is an equilibrium whenever
\[
B^T_s \frac{\partial H}{\partial p}(q,p)=0, \quad B_s\frac{\partial H}{\partial q}(q,p) + B_dRB_d^T\frac{\partial H}{\partial p}(q,p)=0
\]
Premultiplication of the second equation by the row-vector $\frac{\partial^T H}{\partial p}(q,p)$, making use of the first equation, yields $\frac{\partial^T H}{\partial p}(q,p)B_dRB_d^T\frac{\partial H}{\partial p}(q,p)=0$, or equivalently $B_d^T\frac{\partial H}{\partial p}(q,p)=0$. This in turn implies $B_s\frac{\partial H}{\partial q}(q,p) =0$.
\qed

In other words, for $(q,p)$ to be an equilibrium, $\frac{\partial H}{\partial p}(q,p)$ should satisfy the consensus conditions corresponding to the spring-damper graph $\mathcal{G}(\mathcal{V}, \mathcal{E}_s \cup \mathcal{E}_d)$, whereas $\frac{\partial H}{\partial q}(q,p)$ should be in the space of {\it cycles} of the spring graph $\mathcal{G}(\mathcal{V}, \mathcal{E}_s)$ (corresponding to zero net spring forces applied to the masses at the vertices).

Similarly, the {\it Casimirs} (conserved quantities independent of the Hamiltonian $H$) are computed as follows
\begin{proposition}
The Casimir functions are all functions $C(q,p)$ satisfying
\begin{equation}
\frac{\partial C}{\partial p}(q,p) \in \spa \mathds{1}, \quad \frac{\partial C}{\partial q}(q,p) \in \ker B_s
\end{equation}
\end{proposition}
\noindent
\textbf{Proof.}
$C(q,p)$ is a Casimir if 
\[
\begin{bmatrix} \frac{\partial C}{\partial q}(q,p) & \frac{\partial C}{\partial p}(q,p) \end{bmatrix} 
\begin{bmatrix} 0 & B_s^T \\ -B_s & -B_dRB_d^T \end{bmatrix} =0
\]
or equivalently
\[
\frac{\partial^T C}{\partial p}(q,p)B_s = 0, \quad 
\frac{\partial^T C}{\partial q}(q,p)B^T_s + \frac{\partial^T C}{\partial p}(q,p)B_d R B_d^T=0
\]
Postmultiplication of the second equation by $\frac{\partial C}{\partial p}(q,p)$, making use of the first equation, gives the result.
\qed
Therefore all Casimir functions can be expressed as functions of the {\it linear} Casimir functions
\begin{equation}
C(q,p) = \mathds{1}^Tp, \quad C(q,p)=k^Tq,  k \in \ker B_s
\end{equation}
This implies that starting from an arbitrary initial position $(q_0,p_0) \in \Lambda^1 \times \Lambda_0$ the solution of the mass-spring-damper system (\ref{massspringdamper}) will be contained in the affine space
\begin{equation}
\mathcal{A}_{(q_0,p_0)} := \begin{bmatrix} q_0 \\ p_0 \end{bmatrix} + \begin{bmatrix} 0 \\ \ker \mathds{1}^T \end{bmatrix} + \begin{bmatrix} \im B_s^T \\ 0 \end{bmatrix}
\end{equation}
i.e., for all $t$ the difference $q(t)-q_0$ remains in the space of {\it co-cycles} of the spring graph, while $\mathds{1}^Tp(t)=\mathds{1}^Tp_0$.

\subsection{Stability analysis}
Under generic conditions on the Hamiltonian $H(q,p)$, each affine space $\mathcal{A}_{(q_0,p_0)}$ will intersect the set of equilibria $\mathcal{E}$ in a {\it single} point $(q_{\infty},p_{\infty})$, which will qualify as the point of asymptotic convergence starting from $(q_0,p_0)$ (provided there is enough damping present). In order to simplify the statement of the results we will throughout this subsection consider {\it linear} mass-spring systems, corresponding to a quadratic and decoupled Hamiltonian function
\begin{equation}
H(q,p) = \frac{1}{2}q^TKq + \frac{1}{2}p^TGp,
\end{equation}
where $K$ is the positive diagonal matrix of spring constants, and $G$ is the positive diagonal matrix of reciprocals of the masses. It follows that the set of equilibria is given as $\mathcal{E} = \{ (q,p) \in \Lambda^1 \times \Lambda_0 \mid Kq \in \ker B_s, Gp \in  \spa \mathds{1}\}$, while for each $(q_0,p_0)$ there exists a {\it unique} point $(q_{\infty},p_{\infty}) \in \mathcal{E} \cap \mathcal{A}_{(q_0,p_0)}$. In fact, $q_{\infty}$ is given by the spring graph co-cycle/cycle decomposition 
\begin{equation}\label{qinfty}
q_0 = v_0 + q_{\infty}, \quad v_0 \in \im B_s^T \subset \Lambda^1, Kq_{\infty} \in \ker B_s \subset \Lambda_1,
\end{equation}
while $p_{\infty}$ is uniquely determined by\footnote{$Gp_{\infty}=c \mathds{1}$ where the constant $c$ is determined by the initial value vector $p_0$ via the formula $c \sum_{i=1}^N \Pi_{j \neq i}g_j = (\Pi_{i=1}^N g_i) (\sum_{i=1}^N p_{0i})$.}
\begin{equation}\label{pinfty}
Gp_{\infty} \in \spa \mathds{1}, \mathds{1}^Tp_{\infty} = \mathds{1}^Tp_0
\end{equation} 
This leads to the following asymptotic stability theorem. First note that the energy $H(q,p)=\frac{1}{2}q^TKq + \frac{1}{2}p^TGp$ (which obviously is radially unbounded) satisfies
\begin{equation}\label{lyap}
\frac{d}{dt} H(q,p) = - \frac{\partial^T H}{\partial p}(q,p)B_dRB^T_d\frac{\partial H}{\partial p}(q,p)  = -p^TB_dGRB^T_dGp
\leq 0,
\end{equation}
and thus qualifies as a Lyapunov function; showing at least stability. 
\begin{theorem}\label{thm:stability}
Consider a linear mass-spring-damper system with $H(q,p)= \frac{1}{2}q^TKq + \frac{1}{2}p^TGp$, where $K$ and $G$ are diagonal positive matrices. Then for every $(q_0,p_0)$ there exists a unique equilibrium point $(q_{\infty},p_{\infty}) \in \mathcal{E} \cap \mathcal{A}_{(q_0,p_0)}$, determined by (\ref{qinfty}, \ref{pinfty}). Define the spring Laplacian matrix $L_s := B_sKB_s^T$. Then for every $(q_0,p_0)$ the following holds: the trajectory starting from $(q_0,p_0)$ converges asymptotically to $(q_{\infty},p_{\infty})$ if and only if the largest $GL_s $-invariant subspace contained in $\ker B_d^T$ is equal to $\spa \mathds{1}$.
\end{theorem}
{\bf Proof.} By Lasalle's Invariance principle and (\ref{lyap}) the trajectory converges to the largest invariant subspace contained in $\{(q,p) \mid B^T_dGp=0\}$. Differentiation of $B^T_dGp=0$ yields
\[
0=\frac{d}{dt}B^T_dGp= B^T_dG(-B_sKq -B_dRB_d^TGp) = - B^T_dGB_sKq
\]
while further differentiation gives
\[
0= B^T_dGB_sKB_s^TGp = B^T_dGL_sGp
\]
By repeated differentiation one thus concludes that $Gp(t)$ for $t \to \infty$ will converge to the largest $GL_s$-invariant subspace $\mathcal{V}$ contained in $\ker B_d^T$, while $q(t)$ will converge to the subspace $\{q \mid GB_sKq \in \mathcal{V} \}$. Thus if $\mathcal{V}= \spa \mathds{1}$, then $Gp(t) \to Gp_{\infty}$ with $Gp_{\infty} \in \spa \mathds{1}$. Furthermore $q(t) \to q_{\infty}$ with $GB_sKq_{\infty} \in \spa \mathds{1} \subset \ker B_s^T$, and thus $B_s^TGB_sKq_{\infty} =0$ or equivalently $B_sKq_{\infty} = 0$.
\qed
The condition that the largest $GL_s $-invariant subspace contained in $\ker B_d^T$ is equal to $\spa \mathds{1}$ amounts to {\it pervasive damping}: the influence of the dampers spreads through the whole system. 
\begin{remark}
Theorem \ref{thm:stability} is closely related to recent results on partial consensus for double-integrator multi-agent systems \cite{raisch,camlibel}, as will become clear from the discussion in Section \ref{subsec:symmetry}.
\end{remark}
Another feature of the dynamics of the mass-spring-damper system (\ref{massspringdamper}) is its {\it robustness} with regard to constant external (disturbance) forces. Indeed, consider a mass-spring-damper system with boundary masses (see Section \ref{sub:Mass-spring-systems}) and general Hamiltonian $H(q,p)$, subject to {\it constant} forces $\bar{f}_b$
\begin{equation}\label{massspringdamperdist}
\begin{bmatrix} \dot{q} \\ \dot{p} \end{bmatrix} = \begin{bmatrix} 0 & B_s^T \\ -B_s & -B_dRB_d^T \end{bmatrix} \begin{bmatrix} \frac{\partial H}{\partial q}(q,p) \\ \frac{\partial H}{\partial p}(q,p) \end{bmatrix} + \begin{bmatrix} 0 \\ E \end{bmatrix}\bar{f}_b,
\end{equation}
where we {\it assume}\footnote{If the mapping $q \to \frac{\partial H}{\partial q}(q,0)$ is surjective, then there exists for every $\bar{f}_b$ such a $\bar{q}$ if and only if $\im E \subset \im B_s$.} the existence of a $\bar{q}$ such that
\begin{equation}
B_s\frac{\partial H}{\partial q}(\bar{q},0) = E\bar{f}_b
\end{equation}
Then the {\it availability function}
\begin{equation}
\bar{H}(q,p) := H(q,p) - (q-\bar{q})^T\frac{\partial H}{\partial q}(\bar{q},0) - H(\bar{q},0)
\end{equation}
satisfies
\begin{equation}
\frac{d}{dt} \bar{H} (q,p) = - \frac{\partial^T H}{\partial p}(q,p)B_dRB_d^T\frac{\partial H}{\partial p}(q,p) \leq 0
\end{equation}
Specializing to $H(q,p)=\frac{1}{2}q^TKq + \frac{1}{2}p^TGp$, in which case $\bar{H}(q,p) = \frac{1}{2}(q-\bar{q})^TK(q-\bar{q}) + \frac{1}{2}p^TGp$, we obtain the following analog of Theorem \ref{thm:stability}.

\begin{proposition}
Consider a linear mass-spring-damper system (\ref{massspringdamperdist}) with constant external disturbance $\bar{f}_b$ and Hamiltonian $H(q,p)= \frac{1}{2}q^TKq + \frac{1}{2}p^TGp$, where $K$ and $G$ are diagonal positive matrices. and with $\im E \subset \im B_s$. The set of controlled equilibria is given by $\mathcal{\bar{E}} = \{ (q,p) \in \Lambda^1 \times \Lambda_0 \mid B_sKq=E\bar{f}_b, Gp \in \spa \mathds{1}\}$. For every $(q_0,p_0)$ there exists a unique equilibrium point $(\bar{q}_{\infty},p_{\infty}) \in \mathcal{\bar{E}} \cap \mathcal{A}_{(q_0,p_0)}$. Here $p_{\infty}$ is determined by (\ref{pinfty}), while $\bar{q}_{\infty} = \bar{q} + q_{\infty}$, with $\bar{q}$ such that $B_s K\bar{q} = E\bar{f}_b$ and $q_{\infty}$ the unique solution of (\ref{qinfty}) with $q_0$ replaced by $q_0 - \bar{q}$. Furthermore, for each $(q_0,p_0)$ the trajectory starting from $(q_0,p_0)$ converges asymptotically to $(\bar{q}_{\infty},p_{\infty})$ if and only if the largest $GL_s $-invariant subspace contained in $\ker B_d^T$ is equal to $\spa \mathds{1}$.
\end{proposition}
Note that the above proposition has a classical interpretation in terms of the robustness of {\it integral control} with regard to constant disturbances: the springs act as integral controllers which counteract the influence of the unknown external force $\bar{f}_b$ so that the vector of momenta $p$ will still converge to consensus\footnote{The above proposition can be also applied to leader-follower networks (Section \label{sub:clustering}), implying that constant leader inputs can be counteracted by integral control (on top of the standard consensus algorithm).}.

\begin{remark}
The analysis for a mass-damper systems with constant external velocities (\ref{mass-damper}), and equivalently for the leader-follower network (\ref{PHconsensus}) is somewhat different.
Assuming the graph to be {\it connected}, it is well-known \cite{Bollobas98} that for each vector $\bar{e}^b$ there exists a unique
equilibrium vector $\bar{p}$ such that
\begin{equation}
0 = -B_iRB_i^T\bar{p} - B_iRB_b^T\bar{e}^b
\end{equation}
{\it Asymptotic stability} of $\bar{p}$ can then be proved by defining the availability function
$
\bar{H}(p) := H(p) - (p-\bar{p})^T\frac{\partial H}{\partial p}(\bar{p}) - H(\bar{p}) 
$, satisfying
\begin{equation}
\begin{array}{rcl}
\frac{d}{dt}\bar{H}(p) & = & - \begin{bmatrix} (p - \bar{p})^T & (e^b - \bar{e}^b)^T \end{bmatrix} BGB^T 
\begin{bmatrix} p - \bar{p} \\ e^b - \bar{e}^b \end{bmatrix} \\[2mm]
&& + (f_b - \bar{f}_b)^T(e^b - \bar{e}^b)
\end{array}
\end{equation}
where $\bar{f}_b = B_bGB_i^T\bar{p} + B_bGB_b^T\bar{e}^b$ is the output equilibrium value. Since the set
$
\{ p \mid (p - \bar{p})^T B_iGB_i^T (p - \bar{p}) = 0 \}
$
is equal to the single point $\bar{p}$ (because $\ker BGB^T = \spa \mathds{1}$), this shows asymptotic stability of the controlled equilibrium $\bar{p}$ for $e^b = \bar{e}^b$.
\end{remark}

\section{Port-Hamiltonian systems on graphs obtained by symmetry reduction}

In this section we will show how port-Hamiltonian systems on graphs, such as the mass-spring-damper systems, can be alternatively obtained by {\it symmetry reduction} from a \textit{symplectic} formulation, exploiting the invariance of the Hamiltonian function (in particular, of the spring potential energies). 

\subsection{Symmetry reduction from the symplectic formulation}\label{subsec:symmetry}

Let us return to the formulation of a mass-spring system in Section
\ref{sub:Mass-spring-systems}, where the vertices correspond to the masses,
and the edges to the springs in between them. An alternative is to
consider the configuration\footnote{Note that $q_{c}\in\Lambda^{0}$ is defined to be a function $q_{c}:\mathcal{V}\to\mR$, assigning to each vertex its position in $\mathcal{R}$.} vector $q_{c}\in\Lambda^{0}=:Q_{c}$, describing the \textit{positions
of all the masses.} In fact, this is the classical starting point
for Lagrangian mechanics, where we do \textit{not} start with the
\textit{energy variables} $q$ and $p$, but instead we start with
the configuration vector $q_{c}$ and the corresponding velocity vector
$\dot{q}_{c}$. The classical Hamiltonian formulation is then obtained
by \textit{defining} the vector of momenta $p\in\Lambda_{0}=Q_{c}^{*}$
as $p=M\dot{q}_{c}$ (with $M$ the diagonal mass matrix), resulting
in the \textit{symplectic phase space} 
$Q_{c}\times Q_{c}^{*}=\Lambda^{0}\times\Lambda_{0}$.
For ordinary springs the relation between $q_{c}\in\Lambda^{0}$ and the vector
$q\in\Lambda^{1}$ describing the elongations of the springs is given
as $q=B^{T}q_{c}$. Hence in this case the Hamiltonian can be also expressed as a function $H_{c}$ of $(q_{c},p)$ by defining 
\begin{equation}
H_{c}(q_{c},p):=H(B^{T}q_{c},p)\label{Hamiltonianphasespace}
\end{equation}
It follows that the equations of motion of the mass-spring system
(with boundary masses) are given by the canonical Hamiltonian equations
\begin{equation}\label{Hamiltoniansystemphasespace}
\begin{array}{rcl}
\dot{q}_{c} & = & \frac{\partial H_{c}}{\partial p}(q_{c},p)\\[2mm]
\dot{p} & = & -\frac{\partial H_{c}}{\partial q_{c}}(q_{c},p)+ Ef_{b}\\[2mm]
e^{b} & = & E^T \frac{\partial H_{c}}{\partial p}(q_{c},p)
\end{array}
\end{equation}
where, as before, $f_{b}$ are the external forces exerted on the
boundary masses and $e_{b}$ are their velocities.

What is the relation with the port-Hamiltonian formulation given in
Section \ref{sub:Mass-spring-systems} ? It turns out that this relation is precisely
given by the standard procedure of \textit{symmetry reduction} of
a Hamiltonian system\footnote{This relation can be regarded as the discrete, graph-theoretic, version, of the correspondence between the port-Hamiltonian formulation of the Maxwell equations (using the Stokes-Dirac structure) and its symplectic formulation using the vector potential of the magnetic field, cf. \cite{MarsdenRatiu, Vankerschaver}.}. Indeed, since $\mathds{1}^TB=0$ the Hamiltonian
function $H_{c}(q_{c},p)$ given in (\ref{Hamiltonianphasespace})
is \textit{invariant} under the action of the group $\mathfrak{G}=\mR$
acting on the phase space $\Lambda^{0}\times\Lambda_{0} \simeq \mR^{2N}$
by the symplectic group action 
\begin{equation}
(q_{c},p)\mapsto(q_{c}+\alpha\mathds{1},p)\,,\quad\alpha\in\mathfrak{G}=\mR
\end{equation}
From standard reduction theory, see e.g. \cite{MarsdenRatiu, Libermann_marle87} and the references quoted therein, it follows that
we may factor out the configuration space $Q_{c}:=\Lambda^{0}$ to
the \textit{reduced configuration space} 
\begin{equation}
Q:=\Lambda^{0}/\mathfrak{G}
\end{equation}
Let us first assume that the graph is \textit{connected}, in which
case, see e.g. \cite{Bollobas98}, $\ker B^{T}=\spa\mathds{1}$.
Then we have the following identification 
\begin{equation}
Q:=\Lambda^{0}/\mathfrak{G} \simeq B^{T}\Lambda^{0} \subset \Lambda^{1}
\end{equation}
Hence the \textit{reduced state space} of the mass-spring system
is given by $\im B^T \times\Lambda_{0}$, where $\im B^T \subset \Lambda^1$. Furthermore, under the symmetry action the canonical Hamiltonian equations (\ref{Hamiltoniansystemphasespace})
on the symplectic space $\Lambda^{0}\times\Lambda_{0}$ reduce to
the port-Hamiltonian equations (\ref{mass-spring1}) on $\im B^T \times\Lambda_{0} \subset \Lambda^1 \times \Lambda_{0} $
obtained before: 
\begin{equation}\label{mass-spring}
\begin{array}{rcl}
\dot{q} & = & B^{T}\dot{q}_{c}=B^{T}\frac{\partial H_{c}}{\partial p}(q_{c},p)=B^{T}\frac{\partial H}{\partial p}(q,p)\\[2mm]
\dot{p} & = & -\frac{\partial H_{c}}{\partial q_{c}}(q_{c},p)+ Ef_{b}=-B\frac{\partial H}{\partial q}(q,p)+ Ef_{b}\\[2mm]
e^{b} & = & E^T\frac{\partial H}{\partial p}(q,p)
\end{array}\end{equation}
In case the graph is not connected, then the above symmetry reduction
can be performed for each component of the graph (i.e., the symmetry
group is $\mR^{c_{\mathcal{G}}}$, with $c_{\mathcal{G}}$ denoting
the number of components of the graph $\mathcal{G}$), yielding again
the reduced state space\footnote{Note that in fact the subspace $\im B^T \subset\Lambda^1$ is determined by the Casimirs $k^Tq, Bk=0$ in the sense that $\im B^T = \{ q \in \Lambda^0 \mid k^Tq = 0, \mbox{ for all } k \in \ker B \}$. Furthermore, $\im B^T = \Lambda^1$ if and only if the graph does not contain cycles.} $\im B^T \times\Lambda_{0}$.

For a mass-spring-{\it damper} system, although in the standard symmetry reduction framework not considered as a Hamiltonian system, the same reduction procedure can still be applied. A mass-spring-damper system in coordinates $(q_c,p)$ takes the form 
\begin{equation}\label{Hamiltoniansystemphasespace+damper}
\begin{array}{rcl}
\dot{q}_{c} & = & \frac{\partial H_{c}}{\partial p}(q_{c},p)\\[2mm]
\dot{p} & = & -\frac{\partial H_{c}}{\partial q_{c}}(q_{c},p) - B_dRB_d^T\frac{\partial H_{c}}{\partial p}(q_{c},p) + Ef_{b}\\[2mm]
e^{b} & = & E^T \frac{\partial H_{c}}{\partial p}(q_{c},p)
\end{array}
\end{equation}
where $H_c(q_c,p)=H(B_s^Tq,p)$ with $q=B_s^Tq_c$ the spring elongations. Here $B_s$ and $B_d$ denote, as before, the incidence matrices of the spring, respectively, damper graph. Under the same symmetry action as above this reduces to the equations (\ref{massspringdamper}) on the reduced state space $\im B_s^T \times\Lambda_{0}$.

\smallskip

The precise relation between Theorem \ref{thm:stability} and the results obtained in \cite{camlibel, raisch} now becomes clear. Indeed, the double-integrator networks studied in \cite{raisch,camlibel} correspond to linear mass-spring-damper systems with unit masses, unit spring constants and unit damping coefficients, expressed in the position variables $q_c$ and the velocities $\dot{q}_c$, which are equal to the momenta $p$. Thus Theorem \ref{thm:stability} can be seen to be a direct extension of the velocity consensus result expressed in Corollary 10 of \cite{camlibel}. Note on the other hand that, thanks to the systematic use of the port-Hamiltonian structure, the stability treatment given in Section \ref{sec:analysis} is directly extendable to the nonlinear case. Furthermore we obtain the following corollary to Theorem \ref{thm:stability} regarding to 'second-order consensus':
\begin{corollary}
Consider the mass-spring-damper system (\ref{Hamiltoniansystemphasespace+damper}) in coordinates $(q_c,p)$ where we assume the spring graph to be connected. Then for all initial conditions $q_c(t) \to \spa \mathds{1}, p(t) \to \spa \mathds{1}$ if and only the largest $GL_s$-invariant subspace contained in $\ker B_d^T$ is equal to $\spa \mathds{1}$, and moreover $\ker B_s =0$.
\end{corollary}
{\bf Proof.} Compared to Theorem \ref{thm:stability} on 'velocity-consensus' only the condition $\ker B_s =0$ (no cycles in the spring graph) is new. However it follows from the proof of Theorem \ref{thm:stability} that $q(t)=B_s^Tq_c(t) \to q_{\infty}$ with $B_sKq_{\infty}=0$. Thus $q(t)=B_s^Tq_c(t) \to 0$ if and only if $\ker B_s=0$. Finally, by connectedness of the spring graph $q(t)=B_s^Tq_c(t) \to 0$ if and only if $q(t) \to \spa \mathds{1}$.
\qed


\subsection{Further reduction}
It is well-known that symmetry reduction of a Hamiltonian system entails two steps \cite{MarsdenRatiu}. Roughly speaking, the first step, as discussed above, deals with factoring out the state space by the action of the symmetry group, leading to a Hamiltonian system defined with respect to a Poisson structure possessing Casimirs. The second step deals with {\it restricting} the obtained Hamiltonian dynamics to the level sets of these Casimirs; thereby obtaining a {\it reduced symplectic Hamiltonian system}. 

In the case of a mass-spring system (with connected graph) the second step is performed as follows. Note that the dual space $(\im B^T)^*$ can be identified with
\begin{equation}
(\im B^T)^* = \Lambda_1 / (\im B^T)^{\perp} = \Lambda_1 / \ker B \simeq \im B
\end{equation}
Thus the reduced state space can be identified with
\begin{equation}
\im B^T \times \im B \subset \Lambda^1 \times \Lambda_0
\end{equation}
which is again a symplectic space. We leave the extension to the non-connected case and the presence of dampers to the readers.

\section{The Kirchhoff-Dirac structure on graphs and its port-Hamiltonian dynamics}\label{sec:Kirchhoff}

In this section we consider a third canonical graph Dirac structure, which results from {\it constraining} the flows at the internal vertices to {\it zero} (and thus there is no energy-storage or dissipation associated with the vertices for the corresponding port-Hamiltonian system).

\subsection{The Kirchhoff-Dirac structure \label{sub:The-Kirchhoff-Dirac-structure}}

As already alluded to at the end of Section \ref{sub:derived} the \textit{Kirchhoff-Dirac structure} is defined as 
\begin{equation}
\begin{array}{c}
\mathcal{D}_{K}(\mathcal{G}) :=  \{(f_{1},e^{1},f_{b},e^{b})\in\Lambda_{1}\times\Lambda^{1}\times\Lambda_{b}\times\Lambda^{b}\mid\\
\\    B_{i}f_{1}=0,B_{b}f_{1}=f_{b},\,\exists e^{0i}\in\Lambda^{0i}\mbox{ s.t. }e^{1}=-B_{i}^{*}e^{0i}-B_{b}^{*}e^{b}\}
\end{array}\label{KDstructure}
\end{equation}
Note that, in contrast to the flow/effort-continuous graph Dirac structures, the Kirchhoff-Dirac structure only involves the flow
and effort variables of the {\it edge} and {\it boundary} vertex spaces (not of the internal vertex spaces).

\begin{proposition}\label{prop:Kirchhoff} 
$\mathcal{D}_{K}(\mathcal{G})$
is a separable Dirac structure.
\end{proposition} 
\noindent
\textbf{Proof.}
The Kirchhoff-Dirac structure is equal to the composition of the flow-continuous\footnote{Or the composition of the effort-continuous graph Dirac structure with $\{(f_{0},e^{0}) \in\Lambda_{0}\times\Lambda^{0}\mid f_{0}=0\}$.} graph Dirac structure $\mathcal{D}_f(\mathcal{G})$ with the trivial
separable Dirac structure defined as
\[
\{(f_{0i},e^{0i})\in\Lambda_{0i}\times\Lambda^{0i}\mid f_{0i}=0\}
\]
The result then follows from Proposition \ref{compositional}. \qed
\medskip
Port-Hamiltonian systems with respect to the Kirchhoff-Dirac structure are defined completely similar to the case of the flow/effort-continuous graph Dirac structure; the difference being that energy-storing or dissipative relations are now only defined for the flow and effort variables corresponding to the edges.

\subsection{Electrical circuits \label{sub:Electrical-circuits}}

{\it The} example of a port-Hamiltonian system\footnote{The terminology 'port-Hamiltonian' may be confusing in this context, because 'ports' in electrical circuits are usually defined by {\it pairs of terminals}, that is \textit{pairs} of boundary vertices with
external variables being the currents through and the voltages across
an edge corresponding to each such port. See also the discussion in \cite{willemsmodeling, willemsportsterminals,schaftBosgraboek}.} with respect to a Kirchhoff-Dirac structure is an electrical RLC-circuit, with circuit graph $\mathcal{G}$. In this case the elements of $\Lambda_{1}$ and $\Lambda^{1}$ denote
the vectors of currents through, respectively the voltages across,
the edges, and the Kirchhoff-Dirac structure amounts to Kirchhoff's current and voltage laws (whence its name). Furthermore, the effort
variables $e^{0}$ are the \textit{potentials} at
the vertices, while the boundary flows and efforts $f_{b},e^{b}$
are the \textit{boundary currents}, respectively \textit{boundary
potentials} at the boundary vertices (the \textit{terminals}
of the electrical circuit).

On top of Kirchhoff's laws, the dynamics is defined by the energy-storage relations corresponding to either capacitors or inductors, and dissipative relations corresponding to resistors. 
The energy-storing relations for a capacitor at edge $e$ are given by 
\begin{equation}
\dot{Q}_{e}= -I_{e},\quad V_{e}=\frac{dH_{Ce}}{dQ_{e}}(Q_{e})
\end{equation}
with $Q_{e}$ the charge, and $H_{Ce}(Q_{e})$ denoting the electric
energy stored in the capacitor. 
Alternatively, in the case of an inductor one specifies the magnetic
energy $H_{Le}(\Phi_{e})$, where $\Phi_{e}$ is the magnetic
flux linkage, together with the dynamic relations 
\begin{equation}
\dot{\Phi}_{e}=V_{e},\quad - I_{e}= \frac{dH_{Le}}{d\Phi_{e}}(\Phi_{e})
\end{equation}
Finally, a resistor at edge $e$  corresponds to a static relation between the current $I_{e}$
through and the voltage $V_{e}$ across this edge, such
that $V_{e}I_{e}\leq0$. In particular, a linear (ohmic) resistor at edge
$e$ is specified by a relation $V_{e}= -R_{e}I_{e},$ with $R_{e} \geq 0$. 

Alternatively, we can decompose the circuit graph $\mathcal{G}$ as
the interconnection of a graph corresponding to the capacitors, a
graph corresponding to the inductors, and a graph corresponding to
the resistors. For simplicity let us restrict ourselves to the case
of an $LC$-circuit without boundary vertices. Define $\hat{\mathcal{V}}$
as the set of all vertices that are adjacent to at least one capacitor
\textit{as well as} to at least one inductor. Then split the circuit
graph into an open circuit graph $\mathcal{G}^{C}$ corresponding
to the capacitors and an open circuit graph $\mathcal{G}^{L}$ corresponding
to the inductors, both with set of boundary vertices $\hat{\mathcal{V}}$.
Denote the incidence matrices of these two circuit graphs by 
\[
B^{C}:=\begin{bmatrix}B_{i}^{C}\\
B_{b}^{C}\end{bmatrix},\, B^{L}:=\begin{bmatrix}B_{i}^{L}\\
B_{b}^{L}\end{bmatrix}
\]
Assuming for simplicity that all capacitors and inductors are linear
we arrive at the following equations for the $C$-circuit 
\[
\begin{array}{rcl}
B_{b}^{C}\dot{Q} & = & I_{b}^{C},\quad B_{i}^{C}\dot{Q}=0, \\[2mm]
B_{b}^{CT}\psi_{b}^{C} & = & C^{-1}Q-B_{i}^{CT}\psi_{i}^{C}\end{array}
\]
with $Q$ the vector of charges of the capacitors and $C$ the diagonal
matrix with diagonal elements given by the capacitances of the capacitors.
Similarly for the $L$-circuit we obtain the equations 
\[
\begin{array}{rcl}
\dot{\Phi} & = & B_{b}^{LT}\psi_{b}^{L}+B_{i}^{LT}\psi_{i}^{L}\\[2mm]
0 & = & B_{i}^{L}L^{-1}\Phi\\[2mm]
I_{b}^{L} & = & -B_{b}^{L}L^{-1}\Phi\end{array}
\]
with $\Phi$ the vector of fluxes and $L$ the diagonal matrix of
inductances of all the inductors.

The equations of the $LC$-circuit are obtained by imposing the interconnection
constraints $\psi_{b}^{C}=\psi_{b}^{L}=:\psi_{i}$ and $I_{b}^{C}+I_{b}^{L}=0$.
By eliminating the boundary currents $I_{b}^{C},I_{b}^{L}$ one thus
arrives at the differential-algebraic port-Hamiltonian equations 
\[
\begin{array}{l}
\begin{bmatrix}B_{i}^{C} & 0\\
0 & B_{i}^{L}\\
B_{b}^{C} & B_{b}^{L}\end{bmatrix}\begin{bmatrix}-\dot{Q}\\
L^{-1}\Phi\end{bmatrix}=0\\[7mm]
\begin{bmatrix}C^{-1}Q\\
-\dot{\Phi}\end{bmatrix}=\begin{bmatrix}B_{i}^{CT} & 0 & B_{b}^{CT}\\
0 & B_{i}^{LT} & B_{b}^{LT}\end{bmatrix}\begin{bmatrix}\psi_{i}^{C}\\
\psi_{i}^{L}\\
\psi_{i}\end{bmatrix}\end{array}
\]
For a formulation of pure $R, L$ or $C$ circuits, and their weighted Laplacian matrices, we refer to \cite{Schaft_SCL10}.

\subsection{Mass-spring systems with regard to a Lagrangian tree \label{sub:Mass-spring-systems-Lagrangian}}

An alternative port-Hamiltonian formulation of mass-spring(-damper) systems, in terms of the Kirchhoff-Dirac structure, can be given as follows. 
Recall the port-Hamiltonian
formulation on $\Lambda^{1}\times\Lambda_{0}$ with respect to the
effort-continuous graph Dirac structure $\mathcal{D}_{e}(\mathcal{G})$, in which case
the masses correspond to the vertices, and the springs to the edges
of the graph $\mathcal{G}$, which we assume to be connected\footnote{For non-connected graphs $\mathcal{G}$ the same construction can
be done for every connected component.
}.
This graph can be extended to an {\it augmented} graph $\mathcal{G}_{\mathrm{aug}}$
by adding a \textit{ground vertex} $g$ and adding edges from every
vertex $v$ of $\mathcal{G}$ towards this ground vertex. (The augmented
graph is called a \textit{Lagrangian tree}.) Furthermore, by \textit{constraining}
the effort $e_{g}$ at the ground vertex to be zero we can equate
the efforts $e_{v}$ at every vertex $v$ of $\mathcal{G}$ with the
effort $e_{vg}$ at the edge from $v$ to $g$ of the augmented graph
$\mathcal{G}_{\mathrm{aug}}$. In this way we can identify the effort-continuous graph
Dirac structure $\mathcal{D}_{e}(\mathcal{G})$ with the Kirchhoff-Dirac
structure $\mathcal{D}_{K}(\mathcal{G}_{\mathrm{aug}})$ with the additional
constraint $e_{g}=0$. (Note that this is again a separable Dirac
structure since it equals the composition of the Kirchhoff-Dirac structure
$\mathcal{D}_{K}(\mathcal{G}_{\mathrm{aug}})$ with the trivial Dirac structure
$\{(f_{g},e_{g})\mid e_{g}=0\}$.)

In this way, the masses become associated with the edges $e_{vg}$ from
every vertex $v$ to the ground vertex $g$. The interpretation of
the ground vertex $g$ is that it represents the reference point (with
velocity $e_{g}$ being zero). The flow $f_{g}$ at the ground vertex
$g$ equals the total force exerted on a mass located at this reference
point.

\subsection{Properties of the boundary flows and efforts of the Kirchhoff-Dirac structure}

The fact that the internal vertex flows in the definition of the Kirchhoff-Dirac structure are all zero (and consequently no storage or dissipation at the vertices takes place) has a number of specific consequences for the behavior of the boundary flows and efforts (see \cite{willemsportsterminals} for closely related considerations).

Assume (for simplicity of exposition) that $\mathcal{R}=\mathbb{R}$. 
From the definition of the Kirchhoff-Dirac structure and $\mathds{1}^T B=0$ it
follows that 
\begin{equation}
0=\mathds{1}^{T}Bf_{1}=\mathds{1}_{b}^{T}B_{b}f_{1}=-\mathds{1}_{b}^{T}f_{b}\label{constraint}
\end{equation}
with $\mathds{1}_{b}$ denoting the vector with all ones of dimension
equal to the number of boundary vertices. Hence the boundary part
of the Kirchhoff-Dirac structure of an open graph is constrained by
the fact that the boundary flows add up to zero. Dually, we may always
add to the vector of vertex efforts $e^{0}$ the vector $\mathds{1}$
leaving invariant the edge efforts $e^{1}=B^{T}e^{0}$. Hence, to
the vector of boundary efforts $e^{b}$ we may always add the vector
$\mathds{1}_{b}$. \begin{proposition}\label{prop:constraint} Consider
an open graph $\mathcal{G}$ with Kirchhoff-Dirac structure $\mathcal{D}_{K}(\mathcal{G})$.
Then for each $(f_{1},e^{1},f_{b},e^{b})\in\mathcal{D}_{K}(\mathcal{G})$
it holds that 
\[
\mathds{1}_{b}^{T}f_{b}=0 \, ,
\]
while for any constant $c\in\mathbb{R}$ 
\[
(f_{1},e^{1},f_{b},e^{b}+c\mathds{1}_{b})\in\mathcal{D}_{K}(\mathcal{G})
\]
\end{proposition} 
This proposition implies that we may restrict
the dimension of the space of boundary flows and efforts $\Lambda_{b}\times\Lambda^{b}$
of a connected graph by \textit{two}. Indeed, we may define \[
\Lambda_{b\mathrm{red}}:=\{f_{b}\in\Lambda_{b}\mid f_{b}\in\ker\mathds{1}_{b}^{T}\}\]
and its dual space \[
\Lambda_{\mathrm{red}}^{b}:=(\Lambda_{b\mathrm{red}})^{*}=\Lambda^{b}/\im\mathds{1}_{b}\]
It is straightforward to show that the Kirchhoff-Dirac structure
$\mathcal{D}_{K}(\mathcal{G})$ reduces to a linear subspace of the
reduced space $\Lambda_{1}\times\Lambda^{1}\times\Lambda_{b\mathrm{red}}\times\Lambda_{\mathrm{red}}^{b}$,
which is also a Dirac structure. An interpretation of this reduction
is that we may consider one of the boundary vertices, say the first
one, to be the reference vertex, and that we may reduce the vector
of boundary efforts $e^{b}=(e^{b1},\cdots,e^{b\bar{b}})$ to a vector
of voltages $(e^{b2}-e^{b1},\cdots,e^{b\bar{b}}-e^{b1})$. A graph-theoretical
interpretation is that instead of the incidence matrix $B$ we consider
the \textit{restricted} incidence matrix \cite{Bamberg90}.

For a graph $\mathcal{G}$ with more than one connected component
the above holds for each connected component\footnote{The rank of the incidence matrix is equal to the number of vertices minus the number of connected components \cite{Bollobas98}. In fact, each connected component of the graph satisfies the property $\ker B^{T}=\spa\mathds{1}$ with $B$ the (restricted) incidence matrix of this component and
the dimension of $\mathds{1}$ equal to its number of vertices.}. It follows that there
are as many independent constraints on the boundary flows $f_{b}$
as the number of the connected components of the open graph $\mathcal{G}$.
Dually, the space of allowed boundary efforts $e^{b}$ is invariant
under translation by as many independent vectors $\mathds{1}_{b}$
as the number of connected components.

A complementary view on Proposition \ref{prop:constraint} is the
fact that we may \textit{close} an open graph $\mathcal{G}$ to a
closed graph $\bar{\mathcal{G}}$ as follows. Consider first the case
that $\mathcal{G}$ is connected. Then we may add one virtual ground
vertex $v_{g}$, and virtual edges from this virtual vertex to every
boundary vertex $v_{b}\in\mathcal{V}_{e}$, in such a manner that
the Kirchhoff-Dirac structure of this graph $\bar{\mathcal{G}}$ \textit{extends}
the Kirchhoff-Dirac structure of the open graph $\mathcal{G}$. In
fact, to the virtual vertex $v_{g}$ we may associate an arbitrary
potential $e^{0}(v_{g})$ (the ground-potential), and we may rewrite
the externally supplied power $<e^{b}\mid f_{b}>$ as (since by (\ref{constraint})
$\sum_{v_{b}}f_{b}(v_{b})=0$) 
\begin{equation}
<e^{b}\mid f_{b}>=\sum_{v_{b}}(e^{b}(v_{b})-e^{0}(v_{g})f_{b}(v_{b})=\sum_{v_{b}}e^{1b}(v_{b})f_{1b}(v_{b}),
\end{equation}
where $e^{1b}(v_{b}):=e^{b}(v_{b})-e^{0}(v_{g})$ and $f_{1b}(v_{b}):=f_{b}(v_{b})$
denotes the effort across and the flow through the virtual edge towards
the boundary vertex $v_{b}$. It is clear that for every element $(f_{1},e^{1},f_{b},e^{b})\in\mathcal{D}_{K}(\mathcal{G})$
corresponding to the open graph $\mathcal{G}$ there exists $e^{1b}$
such that $(f_{1},e^{1},f_{1b},e^{1b})\in\mathcal{\bar{G}}$ for the
closed graph $\mathcal{\bar{G}}$, and conversely for every $(f_{1},e^{1},f_{1b},e^{1b})\in\mathcal{\bar{G}}$
there exists $e^{b}$ such that $(f_{1},e^{1},f_{b},e^{b})\in\mathcal{D}_{K}(\mathcal{G})$.
This construction is extended to non-connected graphs by adding a ground
vertex to \textit{each} component containing boundary vertices.

\subsection{Physical analogies}

From the above formulation of an RLC-circuit in Section \ref{sub:Electrical-circuits} we conclude that the
structure of the dynamical equations of an inductor is \textit{different}
from that of a capacitor. In order to elucidate this basic difference
we zoom in on the description of an inductor and a capacitor as two-terminal
elements. To this end consider the elementary open graph consisting
of one edge with two boundary vertices $\alpha,\beta$, described
by the incidence matrix $b=\begin{bmatrix}1 &
-1\end{bmatrix}^T$. It follows that an inductor with magnetic energy $H(\Phi)$ is described
by the equations \begin{equation}
\begin{array}{rcl}
\dot{\Phi} & = & b^{T}\begin{bmatrix}\psi_{\alpha}\\
\psi_{\beta}\end{bmatrix}\\[2mm]
\begin{bmatrix}I_{\alpha}\\
I_{\beta}\end{bmatrix} & = & b\frac{dH}{d\Phi}(\Phi),\end{array}\label{inductor}
\end{equation}
whereas a capacitor with electric energy $H(Q)$ is described as 
\begin{equation}
\begin{array}{rcl}
b\dot{Q} & = & \begin{bmatrix}I_{\alpha}\\
I_{\beta}\end{bmatrix}\\[2mm]
\frac{dH}{dQ}(Q) & = & b^{T}\begin{bmatrix}\psi_{\alpha}\\
\psi_{\beta}\end{bmatrix}\end{array}\label{capacitor}
\end{equation}
This difference stems from the fact that the
energy variable $Q$ of a capacitor, as well as the current $I$,
takes values in the linear space $\Lambda_{1}$, while the state variable
$\Phi$ of an inductor, as well as the voltage $V$, takes values
in the {\it dual} space $\Lambda^{1}$. Recalling from Section \ref{sub:massspring} the description of a spring system
\begin{equation}
\begin{array}{rcl}
\dot{q} & = & b^{T}\begin{bmatrix}v_{\alpha}\\
v_{\beta}\end{bmatrix}\\[2mm]
\begin{bmatrix}F_{\alpha}\\
F_{\beta}\end{bmatrix} & = & b\frac{dH}{dq}(q)
\end{array}
\end{equation}
with $q$ the elongation of the spring, and $H(q)$ its potential
energy, we conclude that there is a strict
analogy between a \textit{spring} and an \textit{inductor}\footnote{Thus we favor the so-called \textit{force-current analogy} instead of the \textit{force-voltage analogy}.}. On the other hand, a moving \textit{mass} is \textit{not} a strict
analog of a \textit{capacitor}. Instead, it can be considered to be
the analog of a \textit{grounded} capacitor, while the strict analog of a capacitor (\ref{capacitor}) is the so-called \textit{inerter} \cite{smith} 
\[
b\dot{p} = \begin{bmatrix}F_{\alpha} \\ F_{\beta}\end{bmatrix}, \quad \frac{dH}{dp}(p) =  b^{T}\begin{bmatrix}v_{\alpha}\\
v_{\beta}\end{bmatrix},
\]
where $p$ is the momentum of the inerter and $H(p)$ its kinetic
energy, while $F_{\alpha},F_{\beta}$ and $v_{1},v_{2}$ denote the
forces, respectively, velocities, at the two terminals of the inerter.

\section{Conclusions}

We have laid down a general geometric framework for the description of physical networks dynamics on (non-random) graphs. Starting point are the conservation laws corresponding to the incidence matrix of the graph. This defines three canonical Dirac structures on the combined vertex, edge, and boundary spaces and their duals, where the last one (the Kirchhoff-Dirac structure) corresponds to the absence of energy storage or energy dissipation at the vertices. Relating the internal flows and efforts by either energy-storing or energy-dissipating relations yields various forms of port-Hamiltonian system. We have illustrated the approach on a number of typical physical examples. Other examples that have not been discussed include e.g. supply-chain models and compartmental systems. We have shown how examples from a different origin, such as consensus algorithms, can be formulated and analyzed within the same framework. Furthermore we have shown how classical techniques from Hamiltonian dynamical systems can be exploited for the analysis of the resulting port-Hamiltonian systems.

For clarity of exposition we have only considered the basic building blocks of port-Hamiltonian systems on graphs. Indeed, because the {\it interconnection} of port-Hamiltonian systems again defines a port-Hamiltonian system \cite{DalsmoSIAM99, schaft99, Cervera07}, the framework also covers heterogeneous and multi-scale situations, where several of the constructs considered in the present paper are connected to each other. Moreover, as already indicated in Section \ref{sub:derived} and Remark 3.1, various interesting extensions to dynamical graphs and switching port-Hamiltonian systems on graphs can be made.

The models treated in this paper all correspond to conservation/balance laws within a particular physical domain. Furthermore, the energy-balance of the system components can be seen to {\it result} from the underlying conservation laws and the assumption of integrable constitutive relations for energy-storage. On the other hand, port-based (bond-graph) modeling as originating in the work of Paynter \cite{Paynter61} is aimed at providing a unifying modeling framework for multi-physics systems, by directly {\it starting} from energy-flows between system components from different physical domains. This also results in port-Hamiltonian models as has been amply demonstrated in e.g. \cite{maschkeJFI92, maschkeIEEE_CAS95, schaftAEU95, schaftspringer00, Geoplex09-1}. It is well-known that bond-graph modeling involves an additional abstraction step (e.g., different electrical circuits may lead to the same bond-graph, and, conversely, different bond-graphs may correspond to the same electrical circuit). Furthermore, in the case of electrical circuits port-based modeling starts with a {\it port} description (pairs of terminals), instead of the more basic starting point of {\it terminals} corresponding to conservation laws. Although in most situations the resulting port-Hamiltonian systems are the same this leaves some questions to be answered; see also \cite{willemsportsterminals, schaftBosgraboek}. Another interesting venue for further research \cite{schaftIFAC11} is the precise relation between port-Hamiltonian systems (on graphs) and {\it gradient} dynamical systems; see especially \cite{brayton, smale} for the gradient formulation of RLC circuits.

The identification of the port-Hamiltonian structure, as already crucially used in Section 4 for (stability) analysis, offers important tools for simulation and control. Port-Hamiltonian systems theory has been successful in exploiting the physical structure for control and design purposes, see e.g. \cite{schaftspringer00, OrtegaIEEE_CSM01}, using various forms of passivity-based control, control by interconnection, and tools originating in network synthesis theory. The applications of this control methodology to port-Hamiltonian systems on graphs is an important area for further research. The combination with graph theory, and the inclusion of constraints on the flow and storage variables, is very promising; see already \cite{wei} for preliminary work in this direction.

In a twin paper we will extend the framework from directed graphs to general $k$-complexes. This allows to a give a spatially discretized model of the $2$-D Maxwell equations and of general diffusive systems; see already \cite{schaftCDC08, schaftBosgraboek}. Pertinent questions include the relation of these methods to structure-preserving spatial discretization methods for their description by partial differential equations models. 

\begin{flushleft}
{\bf Acknowledgment} The research of the first author leading to these results has received funding from the  
European Union Seventh Framework Programme [FP7/2007-2013]  under  
grant agreement n257462 HYCON2 Network of Excellence.
\end{flushleft}

\bibliographystyle{plain}

\begin{thebibliography}{10}

\bibitem{Arcak07}
M.~Arcak.
\newblock Passivity as a design tool for group coordination.
\newblock {\em Automatic Control, IEEE Transactions on}, 52(8):1380 --1390,
  2007.

\bibitem{Ball00}
R.S. Ball.
\newblock {\em A {T}reatise on the {T}heory of {S}crews}.
\newblock Cambridge Univerqity Press, Cambridge, UK, 1998.
\newblock First Edition 1900.

\bibitem{Bamberg90}
P.~Bamberg, S.~Sternberg.
\newblock {\em A course in mathematics for students of physics 2}.
\newblock Cambridge University Press, Cambridge, 1990.
\newblock ISBN-0-521-33245.


\bibitem{Bollobas98}
B.~Bollobas.
\newblock {\em Modern Graph Theory}, volume 184 of {\em Graduate Texts in
  Mathematics}.
\newblock Springer, New York, 1998.
\newblock ISBN 0-387-98491-7.

\bibitem{brayton}
R. Brayton, J. Moser.
\newblock A theory of nonlinear networks.
\newblock
{\em Quarterly of {A}pplied {M}athematics}, 22(1,2), 1964, pp. 1--33, 81--104.

%
\bibitem{burger}
M. B\"urger, D. Zelazo, F. Allg\"ower.
\newblock Network clustering: a dynamical systems and saddle-point perspective.
\newblock {\em 50th IEEE Conf. On Decision and Control and European Control Conference (CDC-ECC)}, Orlando, FL, USA, Dec. 12-15, 2011, pp. 7825--7830.

\bibitem{camlibel}
M.K. Camlibel, S. Zhang,
\newblock Partial consensus for heterogenous multi-agent system with double integrator dynamics.
\newblock {\em 4th IFAC Conference on Analysis and Design of Hybrid Systems},
2012, Eindhoven, the Netherlands.

\bibitem{Cervera07}
J.~Cervera, A.J. van~der Schaft, A.~Ba{\~n}os.
\newblock Interconnection of port-{H}amiltonian systems and composition of
  {D}irac structures.
\newblock {\em Automatica}, 43:212--225, 2007.

\bibitem{Courant90}
T.J. Courant.
\newblock Dirac manifolds.
\newblock {\em Trans. American Math. Soc. 319}, pages 631--661, 1990.

\bibitem{DalsmoSIAM99}
M.~Dalsmo, A.J. van~der Schaft.
\newblock On representations and integrability of mathematical structures in
  energy-conserving physical systems.
\newblock {\em SIAM Journal of Control and Optimization}, 37(1):54--91, 1999.

\bibitem{Davies83}
T.~H. Davies.
\newblock Mechanical networks - i: Passivity and redundancy - ii: Formulae for
  the degree of mobility and redundancy - iii: Wrenches on circuit screws.
\newblock {\em Mechanism and Machine Theory}, 18:95--112, 1983.

\bibitem{DePersis}
C.~DePersis, C.S. Kallesoe.
\newblock Pressure regulation in nonlinear hydraulic networks by positive and
  quantized control.
\newblock {\em IEEE Transactions on Control Systems Technology}, 2011.

\bibitem{Dorfman93}
I.~Ya. Dorfman.
\newblock {\em Dirac structures and integrability of nonlinear evolution
  equations}.
\newblock John Wiley, 1993.

\bibitem{Geoplex09-1}
V.~Duindam, A.~Macchelli, S.~Stramigioli, and H.~Bruyninckx eds.
\newblock {\em Modeling and Control of Complex Physical Systems - The
  Port-Hamiltonian Approach}.
\newblock Springer, Sept. 2009.
\newblock ISBN 978-3-642-03195-3.


\bibitem{gerritsen}
K. Gerritsen, A.J. van der Schaft, W.P.M.H. Heemels.
\newblock On switched Hamiltonian systems.
\newblock In {\it Proceedings 15th International Symposium on Mathematical Theory of Networks and Systems (MTNS2002)}, Eds. D.S. Gilliam, J. Rosenthal, South Bend, August 12-16, 2002.

\bibitem{raisch}
D. Goldin, S.A. Attia, J. Raisch
\newblock Consensus for double integrator dynamics in heterogeneous networks.
\newblock 49th IEEE Conf. On Decision and Control (CDC), Atlanta, USA, Dec., 2010, pp. 4504--4510.

\bibitem{Karger_Novak85}
A.~Karger, J.~Novak.
\newblock {\em Space Kinematics and Lie Groups}.
\newblock Gordon and Breach Science Publishers, New York, 1985.
\newblock ISBN 2-88124-023-2.

\bibitem{Libermann_marle87}
P.~Libermann, C.-M. Marle.
\newblock {\em Symplectic {G}eometry and {A}nalytical {M}echanics}.
\newblock D. Reidel Publishing Company, Dordrecht, Holland, 1987.
\newblock ISBN 90-277-2438-5.


\bibitem{MarsdenRatiu}
J.E. Marsden, T.S. Ratiu.
\newblock {\em Introduction to {M}echanics and {S}ymmetry: a basic exposition
  to Classical Mechanics}.
\newblock London Mathematical Society Lecture Notes Series. Springer, New York,
  2nd edition, 1999.

\bibitem{maschkeLondon97_tot}
B.M. Maschke, A.J. van~der Schaft.
\newblock {\em Modelling and Control of Mechanical Systems}, chapter
  Interconnected Mechanical systems. Part 1 and 2, pages 1--30.
\newblock Imperial College Press, London, 1997.
\newblock ISBN 1-86094-058-7.

\bibitem{maschkeJFI92}
B.M. Maschke, A.J. van~der Schaft, P.C. Breedveld.
\newblock An intrinsic {H}amiltonian formulation of network dynamics:
  Non-standard {P}oisson structures and gyrators.
\newblock {\em Journal of the Franklin institute}, 329(5):923--966, 1992.

\bibitem{maschkeIEEE_CAS95}
B.M. Maschke, A.J. van~der Schaft, P.C. Breedveld.
\newblock An intrinsic {H}amiltonian formulation of the dynamics of
  {LC}-circuits.
\newblock {\em IEEE Trans. Circuits and Systems I:Fundamental Theory and
  Applications}, 42(2):73--82, 1995.
  
\bibitem{twin} B. Maschke, A.J. van der Schaft, in preparation.

\bibitem{murray94}
R.~M. Murray, X.~Li, S.S. Sastry.
\newblock {\em A Mathematical Introduction to Robotic Manipulation}.
\newblock CRC press, March 1994.
\newblock ISBN 0-8493-7981-4.

\bibitem{Olfati-SaberProcIEEE07}
R.~Olfati-Saber, J.A. Fax, R.M. Murray.
\newblock Consensus and cooperation in networked multi-agent systems.
\newblock {\em Proceedings of the IEEE}, 95(1):215 --233, 2007.
%

\bibitem{OrtegaIEEE_CSM01}
R.~Ortega, A.J. van~der Schaft, I.~Mareels, B.~Maschke.
\newblock Putting energy back in control.
\newblock {\em IEEE Control Systems Magazine}, 21(2):18-- 32, April 2001.

\bibitem{Paynter61}
H.~M. Paynter.
\newblock {\em Analysis and design of engineering systems}.
\newblock M.I.T. Press, Cambridge, Mass., 1961.

\bibitem{Rahmani09}
A.~Rahmani, M.~Ji, M.~Mesbah, M.~Egerstedt.
\newblock Controllability of multi-agent systems from a graph-theoretic
  perspective.
\newblock {\em SIAM J. Control Optim.}, 48(1):162--186, 2009.

\bibitem{bookhydraulic}
J.A. Roberson, C.T. Crowe.
\newblock {\em Engineering fluid mechanics}.
\newblock Houghton Mifflin Company, 1993.

\bibitem{Selig96}
J.M. Selig.
\newblock {\em {G}eometric {M}ethods in {R}obotics}.
\newblock Monographs in {C}omputer {S}ciences. Springer Verlag, 1996.
\newblock ISBN 0-387-94728-0.


\bibitem{schaft99}
A.J. van~der Schaft.
\newblock {\em The Mathematics of Systems and Control, From Intelligent Control
  to Behavioral Systems}, chapter Interconnection and {G}eometry, pages
  203--218.
\newblock eds. J.W. Polderman, H.L. Trentelman. Un. Groningen, Groningen, The
  Netherlands, 1999.

\bibitem{schaftspringer00}
A.J. van~der Schaft.
\newblock {\em $L_2$-{G}ain and {P}assivity {T}echniques in {N}onlinear
  Control}.
\newblock Springer Communications and Control Engineering series.
  Springer-Verlag, London, 2nd revised and enlarged edition, 2000.
\newblock first edition Lect. Notes in Control and Inf. Sciences, vol. 218,
  Springer-Verlag, Berlin, 1996.

\bibitem{Schaft_SCL10}
A.J. van~der Schaft.
\newblock Characterization and partial synthesis of the behavior of resistive
  circuits at their terminals.
\newblock {\em Systems Control Letters}, 59(7):423 -- 428, 2010.

\bibitem{schaftIFAC11}
A.J. van~der Schaft.
\newblock On the relation between port-Hamiltonian and gradient systems.
\newblock In {\em Proc. IFAC World Congress}, Milan, Italy, 2011.

\bibitem{schaftcamlibel}
A.J. van der Schaft, M.K. Camlibel.
\newblock A state transfer principle for switching port-Hamiltonian systems.
\newblock In {\it Proc. 48th IEEE Conf. on Decision and Control}, Shanghai, China, December 16-18, 2009, pp. 45--50.

\bibitem{schaftAEU95}
A.J. van~der Schaft, B.M. Maschke.
\newblock The {H}amiltonian formulation of energy conserving physical systems
  with external ports.
\newblock {\em Archiv f{\"u}r Elektronik und {\"U}bertragungstechnik},
  49(5/6):362--371, 1995.

\bibitem{schaftGeomPhys02}
A.J. van~der Schaft and B.M. Maschke.
\newblock Hamiltonian formulation of distributed parameter systems with
  boundary energy flow.
\newblock {\em J. of Geometry and Physics}, 42:166--174, 2002.

\bibitem{schaftCDC08}
A.J. van~der Schaft, B.M. Maschke.
\newblock Conservation laws on higher-dimensional networks.
\newblock In {\em Proc. 47th IEEE Conf. on Decision and Control}, Cancun,
  Mexico, 2008.

\bibitem{schaftBosgraboek}
A.J. van~der Schaft, B.M. Maschke.
\newblock {\em Model-Based Control: Bridging Rigorous Theory and Advanced
  Technology, P.M.J. Van den Hof, C. Scherer, P.S.C. Heuberger, eds.}, chapter
  Conservation laws and lumped system dynamics, pages 31--48.
\newblock Springer, Berlin-Heidelberg, 2009.

\bibitem{schaftNECSYS10}
A.J. van~der Schaft, B.M. Maschke.
\newblock Port-{H}amiltonian dynamics on graphs: Consensus and coordination
  control algorithms.
\newblock In {\em Proc. 2nd IFAC Workshop on Distributed Estimation and Control
  in Networked Systems}, Centre de Congres de L Imp\'erial Palace, Annecy,
  France, 2010.
  
\bibitem{wei} A.J. van der Schaft, J. Wei.
\newblock A Hamiltonian perspective on the control of dynamical distribution
networks.
\newblock to appear in Proceedings 4th IFAC Workshop on Lagrangian and Hamiltonian Methods in Non Linear Control (LHMNLC2012), Bertinoro, Italy, August 29--31, 2012.

  \bibitem{smale}
S. Smale.
\newblock On the Mathematical Foundations of Electrical Circuit Theory.
\newblock {\em J. of Differential Geometry}, 7, pp. 193--210, 1972.

\bibitem{smith}
M.C. Smith, 
\newblock Synthesis of Mechanical Networks: the Inerter.
\newblock {\it IEEE Trans. Autom. Control}, 47, no. 10, pp. 1648--1662, 2002.


\bibitem{Vankerschaver}
J.~Vankerschaver, H.~Yoshimura, J.~E. Marsden.
\newblock Stokes-{D}irac structures through reduction of infinite-dimensional
  {D}irac structures.
\newblock In {\em Proc. 49th IEEE Conference on Decision and Control}, Atlanta,
  USA, December 2010.
  
\bibitem{willemsmodeling}
J.C. Willems.
\newblock The behavioral approach to open and interconnected systems.
\newblock {\em IEEE Control Systems Magazine}, 27:46--99, 2007.

\bibitem{willemsportsterminals}
J.C. Willems.
\newblock Terminals and ports.
\newblock {\em IEEE Circuits and Systems Magazine}, 10(4):8--16, December 2010.

\bibitem{zelazo}
D. Zelazo, M. Mesbahi.
\newblock Edge agreement: graph-theoretic performance bounds and passivity analysis.
\newblock {\em IEEE Trans. Automatic Control}, vol. 56(3), pp. 544--555, 2011.

\end{thebibliography}

\end{document}